\documentclass[a4paper,11pt,reqno]{amsart}
\usepackage{amsmath}
\usepackage{amssymb,amsfonts,amsthm,amscd}
\usepackage{hyperref}
\numberwithin{equation}{section}

\newtheorem{theorem}{Theorem}[section]
\newtheorem{corollary}[theorem]{Corollary}
\newtheorem{lemma}[theorem]{Lemma}
\newtheorem{proposition}[theorem]{Proposition}
\theoremstyle{definition}
\newtheorem{defn}[theorem]{Definition}

\theoremstyle{remark}
\newtheorem{remark}[theorem]{Remark}

\newcommand{\nc}{\newcommand}
\newcommand{\be}{\begin{equation}}
\newcommand{\ee}{\end{equation}}

\newcommand{\bc}{\begin{center}}
\newcommand{\ec}{\end{center}}
 \nc{\bth}{\begin{theorem}}
 \nc{\bpr}{\begin{proposition}}
 \nc{\epr}{\end{proposition}}
 \nc{\ble}{\begin{lemma}}
 \nc{\ele}{\end{lemma}}
 \nc{\bco}{\begin{corollary}}
 \nc{\eco}{\end{corollary}}
 \nc{\bre}{\begin{remark}}
 \nc{\ere}{\end{remark}}
   \nc{\na}{\nabla}
   \nc{\al}{\alpha}
   \nc{\bet}{\beta}
\begin{document}
\title[
Pseudo-Riemannian submersions from pseudo-hyperbolic
spaces]{Classification of Pseudo-Riemannian submersions with totally
geodesic fibres from pseudo-hyperbolic spaces}
 \subjclass[2010]{53C50}
 \keywords{Pseudo-Riemannian submersions, Pseudo-hyperbolic spaces,
 Special Osserman pseudo-Riemannian manifolds.}
\date{}
\author{Gabriel B\u adi\c toiu}
\address{``Simion Stoilow'' Institute of Mathematics
 of the Romanian Academy, Research Unit 4, P.O. Box 1-764,
 014700 Bucharest, Romania.}
 \email{Gabriel.Baditoiu@imar.ro}
\thanks{The author was supported by a grant of the Romanian National Authority for Scientific
Research, CNCS - UEFISCDI, project number PN-II-ID-PCE-2011-3-0362.}

\begin{abstract}
We classify pseudo-Riemannian submersions with connected totally
geodesic fibres from a real pseudo-hyperbolic space onto a
pseudo-Riemannian manifold. Also, we obtain the classification of
the pseudo-Riemannian submersions with (para-)complex connected
totally geodesic fibres from a (para-)complex pseudo-hyperbolic
space onto a pseudo-Riemannian manifold.
\end{abstract}
\maketitle

\section{Introduction and the main theorem}\label{s1}
\noindent Riemannian submersions, introduced by O'Neill \cite{one} and Gray
\cite{gra}, have been used by many authors to construct new specific
Riemannian metrics, like Einstein or positively curved ones
\cite{bes,gromoll-walschap}, and to study various geometric
structures of Riemannian manifolds \cite{fal-pas-ian}. In this
paper, we show that the pseudo-Riemannian submersions with
connected, totally geodesic fibres from a pseudo-hyperbolic onto a
pseudo-Riemannian manifold are equivalent to the Hopf ones, see
below. First, we give a short review of well-known classification
results of Riemannian submersions and of their nice applications in
Riemannian geometry and then we discuss the pseudo-Riemannian case.

 In early work,
Escobales \cite{esc, esco} and Ranjan \cite{ran} classified
Riemannian submersions with connected totally geodesic fibres from a
sphere, and with complex connected totally geodesic fibres from a
complex projective space. Using a topological argument, Ucci
\cite{ucc} showed that there are no Riemannian submersions with
fibres $\mathbb CP^3$ from the complex projective space $\mathbb
CP^7$ onto $S^8(4)$, and with fibres $\mathbb HP^1$ from the
quaternionic projective space $\mathbb HP^3$ onto $S^8(4)$. A major
advance obtained by Gromoll and Grove in \cite{gro} is that, up to
equivalence, the only Riemannian submersions from spheres with
connected fibres are the Hopf fibrations, except possibly for
fibrations of the $15$-sphere by homotopy $7$-spheres. This
classification was invoked in the proofs of the Diameter Rigidity
Theorem in Gromoll and Grove \cite{grom} and of the Radius Rigidity
Theorem in Wilhelm \cite{wil}. Using Morse theory, Wilking
\cite{wilk}  ruled out the Gromoll and Grove unsettled case by
showing that any Riemannian submersion $\pi:S^{15}\to B^8$ is
equivalent to a Riemannian submersion with totally geodesic fibres,
which by Escobales' classification must be equivalent to a Hopf
Riemannian submersion. A nice consequence of this classification is
the improved version of the Diameter Rigidity Theorem due to Wilking
\cite{wilk}.

In the pseudo-Riemannian set-up, the pioneering work is due to
 Magid  \cite{mag}, who proved that the pseudo-Riemannian submersions with
connected totally geodesic fibres from an anti-de Sitter space onto
 a Riemannian manifold are equivalent to the Hopf
pseudo-Riemannian submersions $H^{2m+1}_1\to\mathbb CH^m$.
Generalizing Magid's result, Stere Ianu\c s and I showed that any
pseudo-Riemannian submersion with connected totally geodesic fibres
from a pseudo-hyperbolic space onto a Riemannian manifold is
equivalent to one of the Hopf pseudo-Riemannian submersions:
$H^{2m+1}_1\to\mathbb CH^m$, $H^{4m+3}_3\to\mathbb HH^m$ or
$H^{15}_7\to\mathbb H^8(-4)$, and as a consequence we classified the
pseudo-Riemannian submersions with connected complex totally
geodesic fibres from a complex pseudo-hyperbolic space onto a
Riemannian manifold (see \cite{dga}).
 In \cite{tohoku}, I extended these results to the case of a pseudo-Riemannian
 base under the assumption that either
 (i) the base space is isotropic or
 (ii) the dimension of fibres is less than or equal
 to $3$, and the metrics induced on the fibres are negative
 definite.
 I also proved that condition (ii) implies (i) (see \cite{tohoku}).
 In this paper, we drop
 these assumptions and we prove the
 following main result.

\begin{theorem}\label{t:main}
    Let $\pi:H^{a}_{l}\to B$ be a pseudo-Riemannian submersion with connected totally
    geodesic fibres from a real pseudo-hyperbolic space
    $H^{a}_{l}$ of curvature $-1$ onto a pseudo-Riemannian manifold.
   Then $\pi$ is equivalent to one of the following
    Hopf pseudo-Riemannian submersions$:$
\begin{itemize}
    \item[$(a)$]  $\pi_{\mathbb C}:H^{2m+1}_{2t+1}\to\mathbb CH^m_t,\ \ 0\leq t\leq m$,
    \ \ \    \ \ \ $(b)$   $\pi_{\mathbb A}:H^{2m+1}_{m}\to\mathbb AP^m$,
    \item[$(c)$]  $\pi_{\mathbb H}:H^{4m+3}_{4t+3}\to\mathbb HH^m_t,\ \ 0\leq t\leq m$,
   \ \ \     \ \ \ $(d)$   $\pi_{\mathbb B}:H^{4m+3}_{2m+1}\to\mathbb BP^m$,
    \item[$(e)$]  $\pi^1_{\mathbb O}:H^{15}_{15} \to H^{8}_{8}(-4) $,
    \ \ \ \ \  \ \ \
     $(f)$   $\pi_{\mathbb O'}:H^{15}_{7}\to H^8_{4}(-4)$, \ \ \  \ \ \  \ \
     $(g)$   $\pi^2_{\mathbb O}:H^{15}_{7}\to H^8(-4)$.
\end{itemize}
 where $\mathbb CH^m_t$, $\mathbb HH^m_t$ are the indefinite complex
 and quaternionic pseudo-hyperbolic spaces of holomorphic,
 respectively, quaternionic curvature $-4$;
  $\mathbb AP^m$ is the
 para-complex projective space of real dimension  $2m$,  signature
 $(m,m)$, and of para-holomorphic  curvature  $-4$;
  $\mathbb BP^m $ is the  para-quaternionic
 projective space of real dimension  $4m$,   signature  $(2m,2m)$,
 and of para-quaternionic curvature $-4$.
\end{theorem}

The plan of the paper can be summarized as follows. Section \ref{s2}
presents some known definitions and results in the theory of
pseudo-Riemannian submersions. In \S\ref{s3}, we exhibit the
construction of the Hopf pseudo-Riemannian submersions from
pseudo-hyperbolic spaces, which ensures the existence of at least
one pseudo-Riemannian submersion in each class (a)--(g) of
Theorem \ref{t:main}. In \S\ref{s4}, we see that the  base space $B$
is isometric to either a pseudo-hyperbolic space or a  complete,
simply connected, special Osserman pseudo-Riemannian manifold, which
was classified in \cite{special-osserman-classification}.  To
exclude the Cayley planes of octonions, and of para-octonions from
the list of possible base spaces, we prove that the curvature tensor
of $B$ has a Clifford structure.  For the remaining cases, we
establish that the dimension and the index of the total space are,
in fact, those claimed in Theorem \ref{t:main}. This reduces the
equivalence problem of two pseudo-Riemannian submersions to the one
of the same base space, which  we resolve in \S\ref{s5}. Section
\ref{s6} features consequences of Theorem \ref{t:main}: (a) the
classification of the pseudo-Riemannian submersions with totally
geodesic fibres from complex pseudo-hyperbolic spaces  or from
para-complex projective spaces under the assumption that the fibres are,
respectively, complex or para-complex submanifolds and (b) the
non-existence of the pseudo-Riemannian submersions with quaternionic
or para-quaternionic fibres from $\mathbb HH^m_t$ or $\mathbb
BP^m$, respectively.

\section{Preliminaries}\label{prelim}\label{s2}
In this section we recall several notions and results that will be
used throughout the paper.
\begin{defn}
    A smooth surjective submersion
         $\pi :(M,g)\to (B,g')$
    between two pseudo-Riemannian manifolds is said to be a
    {\it pseudo-Riemannian submersion} (see \cite{onei})
    when $\pi_*$ preserves scalar products of
    vectors normal to fibres and when the metric induced
    on every fibre $F_b=\pi^{-1}(b)$, where $b\in B$, is
    non-degenerate.
\end{defn}
   The vectors tangent to fibres are called
    vertical and those normal to fibres
    are called horizontal.
    We denote  the vertical distribution by $\mathcal V$
    and
    the horizontal distribution  by $\mathcal H$.
    The geometry of  pseudo-Riemannian submersions is characterized in
    terms of the
    O'Neill tensors $T$, $A$ (see \cite{one,onei})
    defined for every vector fields  $E$, $F$ on $M$ by
\begin{eqnarray}
        A_EF =h \na _{hE}{vF} + v \na _{hE}{hF},\ \ \ \ \ \ \ \
        T_EF =h \na _{vE}{vF} + v \na _{vE}{hF},
\end{eqnarray}
    where $\na $ is the Levi-Civita connection of $g$, and $v$ and $h$ denote the orthogonal
    projections on $\mathcal V$  and $\mathcal H$, respectively.
    We  assume that the fibres are totally
    geodesic, which is equivalent to $T_EF=0$ for every $E,F$.
    The O'Neill tensor $A$ is alternating, i.e. $A_XY=-A_YX$ for any
    horizontal vectors $X,Y$, and skew-symmetric with respect to $g$,
    i.e. $g(A_EF,G)=-g(F,A_EG)$ for every vector fields $E$, $F$, $G$
    (see \cite{bes,fal-pas-ian,one,onei}).
    Throughout the paper, $X,Y,Z,Z'$ will always be horizontal
    vector fields, while $U,V,W,W'$ will be vertical vector fields.
    We assume that $\dim M>\dim B$ and that $M$ is connected.

   We denote by  $R$, $R'$ and  $\hat R$ the Riemann curvature tensors
of $M$, $B$ and of a fibre $F_b$, respectively. We choose the
convention for the curvature tensor
$R(E,F)=\na_E\na_F-\na_F\na_E-\na_{[E,F]}$.
By $R'(X,Y)Z$ we shall also
denote  the horizontal lift of $R'(\pi_*X,\pi_*Y)\pi_*Z$.
    The structure equations of a pseudo-Riemannian
    submersion, usually called the O'Neill equations, are
    stated next in a totally geodesic fibre set-up.
\begin{proposition}[\cite{bes,fal-pas-ian,gra,one}]\label{p:oneill}
If $\pi:M\to B$ is a pseudo-Riemannian submersion with totally
geodesic fibres, then

$(a)$
$R(X,Y,Z,Z')=R'(X,Y,Z,Z')-2g(A_XY,A_ZZ')+g(A_YZ,A_XZ')-g(A_XZ,A_YZ'),$

$(b)$ $R(X,Y,Z,U)=g((\nabla_ZA)_XY,U)$,

$(c)$ $R(X,U,Y,V)=g((\nabla_UA)_XY,V)+g(A_XU,A_YV)$,

$(d)$  $R(U,V,W,W')=\hat R(U,V,W,W')$, and $(e)$  $R(U,V,W,X)=0$.
\end{proposition}
\begin{corollary}\label{p:oneill2}
If $\pi:M\to B$ is a pseudo-Riemannian submersion with totally
geodesic fibres, then

 $(a)$  $R(X,Y,X,Y)=R'(X,Y,X,Y)-3g(A_XY,A_XY)$,

 $(b)$  $R(X,U,X,U)=g(A_XU,A_XU)$.
\end{corollary}

\begin{defn}
   A vector field $X$ on $M$ is said to be {\it basic}
    if $X$ is horizontal and $\pi$-related to
    a vector field $X'$ on $B$.
    A vector field $X$ along the fibre
    $\pi^{-1}(b)$, $b\in B$ is said to be {\it basic along}
    $\pi^{-1}(b)$
    if  $X$ is horizontal and $\pi_{*p}X(p)=\pi_{*q}X(q)$
    for every $p$, $q\in \pi^{-1}(b)$.
\end{defn}
We note that each vector field $X'$ on $B$ has a unique horizontal
lift $X$ to $M$ which is basic. For a vertical vector field $V$ and
a  basic vector field $X$ we have $h\na_VX=A_XV$ (see \cite{one}).

\begin{defn}
     Two pseudo-Riemannian submersions $\pi,\pi':(M,g)\to (B,g')$
     are said to be {\it equivalent} if there exists an
     isometry $f$ of $M$ that induces an isometry $\tilde f$ of $B$ so that
     $\pi'\circ f=\tilde f\circ\pi$.
\end{defn}

\section{The construction of the Hopf pseudo-Riemannian
submersions}\label{s3}

In this section, we exhibit the constructions of the real, complex,
quaternionic pseudo-hyperbolic spaces, of the para-complex and
para-quaternionic projective spaces and of the Hopf
pseudo-Riemannian submersions from the real pseudo-hyperbolic
spaces.
\begin{defn}
   Let $\langle\cdot,\cdot\rangle_{\mathbb R^{m+1}_{t+1}}$
   be the inner product of signature $(m-t,t+1)$ on $\mathbb R^{m+1}$
   given by
\begin{eqnarray}
\langle x,y\rangle_{\mathbb
R^{m+1}_{t+1}}=-\sum\limits_{i=0}^tx_iy_i+\sum\limits_{i=t+1}^mx_iy_i
\end{eqnarray}
  for $x=(x_0,\cdot\cdot\cdot,x_m),y=(y_0,\cdot\cdot\cdot,y_m)\in\mathbb R^{m+1}$.
  For any $c<0$ and any positive integer $t$, let
  $H^m_t(c)=\{x\in\mathbb R^{m+1}\ |\ \langle x,x\rangle_{\mathbb R^{m+1}_{t+1}}=1/c\}$
   be the pseudo-Riemannian submanifold of
    $$\mathbb R_{t+1}^{m+1}=(\mathbb R^{m+1},
  ds^2=-dx_0\otimes dx_0-\cdot\cdot\cdot -dx_t\otimes dx_t+
        dx_{t+1}\otimes dx_{t+1}+\cdot\cdot\cdot +dx_m\otimes dx_m).$$
  The space $H^m_t(c)$ is called the $m$-dimensional
 {\it pseudo-hyperbolic space} of index $t$. The hyperbolic space $H^m(c)$ is
 the hypersurface
  $\{x=(x_0,x_1,\cdots,x_m)\in\mathbb R^{m+1}\ |\ x_0>0,\ \langle x,x\rangle_{\mathbb R^{m+1}_{1}}=1/c\}$
  endowed with the  metric induced from $\mathbb R^{m+1}_1$.
\end{defn}

The space $H^m_t(c)$ has constant sectional curvature $c$,
and  we shall define simply  $H^m_t=H^m_t(-1)$.

Throughout the paper, we use the notation: $\mathbb H$ for the
field of quaternions; $\mathbb A$ and $\mathbb B$ for the algebras of
para-complex and para-quaternionic numbers, respectively; $\mathbb
O$ for the algebra of octonions (Cayley numbers)  and $\mathbb O'$
for that of para-octonions \cite{held-stavrov-vankoten} (split
octonions). For $F\in\{\mathbb C, \mathbb A, \mathbb H, \mathbb B,
\mathbb O, \mathbb O'\},$ and for
 $z\in F$, we denote by $\bar{z}$ the conjugate of $z$ in $F$ and, as
 usual,
  $|z|^2_F=\bar{z}z=z\bar z\in\mathbb R$.

\subsection{The indefinite Hopf pseudo-Riemannian
submersions}\label{s:3.1} When $K\in\{\mathbb C, \mathbb H\}$, let
 $\langle\cdot,\cdot\rangle_{K^{m+1}_{t+1}}$
 be the inner product on
 $K^{m+1}$
   given by
\begin{eqnarray}
\langle z,w\rangle_{K^{m+1}_{t+1}}
  =\mathrm{Re}(-\sum\limits_{i=0}^tz_i\bar w_i+\sum\limits_{i=t+1}^mz_i\bar
  w_i),
\end{eqnarray}
where $z=(z_0,\cdot\cdot\cdot,z_m),y=(w_0,\cdot\cdot\cdot,w_m)\in
K^{m+1}$. We set $d=\dim_{\mathbb R}K$ and assume $c<0$. We simply
note that
 $S^{d-1}=\{z\in K\ |\
z\bar z=1\}$,
 and
\begin{eqnarray}
H^{d(m+1)-1}_{d(t+1)-1}(c/4)=\{z\in K^{m+1}\ |\ \langle
z,z\rangle_{K^{m+1}_{t+1}}=4/c\}.
\end{eqnarray}
 The restriction of the projection
\begin{eqnarray}
\ \{z\in K^{n+1} \ |\ \langle z,z\rangle_{K^{m+1}_{t+1}}<0\}
    \to
    \{z\in K^{n+1} \ |\ \langle z,z\rangle_{K^{m+1}_{t+1}}<0\}/K^*,\
 z\mapsto zK^*
\end{eqnarray}
 to
$H^{d(m+1)-1}_{d(t+1)-1}(c/4)$ is a submersion
\begin{eqnarray}
 \pi_K:
 H^{d(m+1)-1}_{d(t+1)-1}(c/4)\to
 KH^m_t(c)=H^{d(m+1)-1}_{d(t+1)-1}(c/4)/S^{d-1},\
 z\mapsto zS^{d-1},
\end{eqnarray}
 called the indefinite Hopf fibration of
$H^{d(m+1)-1}_{d(t+1)-1}(c/4)$. There is a unique pseudo-Riemannian
metric on $KH^{m}_t(c)$ such that $\pi_K:
H^{d(m+1)-1}_{d(t+1)-1}(c/4)\to KH^{m}_t(c)$ is a pseudo-Riemannian
submersion with totally geodesic fibres. We shall simply
 define $KH^{m}_t=KH^{m}_t(-4)$. For $c=-4$, and for $K=\mathbb C$ and  $K=\mathbb
H$,  respectively, the Hopf pseudo-Riemannian submersions are:
\begin{itemize}
\item[(a)]  $\pi_{\mathbb C}:H^{2m+1}_{2t+1}\to\mathbb CH^{m}_t$ with the fibres
isometric to
  $H^1_1=(S^1,-g_{S^1})$, and
\item[(b)] $\pi_{\mathbb H}:H^{4m+3}_{4t+3}\to\mathbb HH^{m}_t$ with the fibres
$H^3_3=(S^3,-g_{S^3})$.
\end{itemize}
A nice reference for the construction of
 $\pi_{\mathbb C}$ is
 \cite{barros-romero}.
 Note that $\mathbb CH^{m}_t$ has holomorphic sectional
 curvature $-4$ (see \cite{barros-romero}),
 and that $\mathbb HH^{m}_t$ has quaternionic sectional
 curvature $-4$.

\subsection{The para-Hopf pseudo-Riemannian
submersions}\label{s:3.2} There are several models of  para-complex
and of para-quaternionic projective spaces
\cite{erdem,gadea-amilibia,cruceanu,blazic}. Following
\cite{erdem,gadea-amilibia}, we  present a para-complex model of a
para-complex projective space, $\mathbb AP^m$, which is simply
connected for $m\geq 2$, and  a simply connected para-quaternionic
model for the para-quaternionic projective space, $\mathbb BP^m$; see
\cite{blazic}.

For $D\in\{\mathbb A, \mathbb B\}$, let $d=\dim_{\mathbb R} D$. We
consider the inner product of signature $(\frac{(m+1)d}{2},$ $
\frac{(m+1)d}{2})$ on $D^{m+1}$ given by
\begin{eqnarray}
\langle z,w\rangle
  =\mathrm{Re}(\sum\limits_{i=0}^mz_i\bar w_i)
\end{eqnarray}
  for $z=(z_0,\cdots,z_m),y=(w_0,\cdots,w_m)\in D^{m+1}.$
Identifying $D^{m+1}=\mathbb R^{d(m+1)}_{d(m+1)/2}$, via
 $(z_0,\cdots,z_m)\simeq(z_0^1,\cdots,z_m^1,\cdots,z_0^d,\cdots,z_m^d)$,
 where $z_i=(z_i^1,\cdots,z_i^d)$, $0\leq i\leq m$,
we simply have $\langle z,w\rangle=-\langle z,w\rangle_{\mathbb
R^{d(m+1)}_{d(m+1)/2}}$, for any $z,w$. In particular, we can write
$H^{2m+1}_{m}=\{z\in\mathbb A^{m+1}\ |\
 \langle z,z\rangle=1 \}$ and
$H^{4m+3}_{2m+1}=\{z\in\mathbb B^{m+1}\ |\
 \langle z,z\rangle=1 \}$.

We set $\mathbb A_0^{m+1} =
 \{z\in\mathbb A^{m+1}\ |\ \langle z,z\rangle>0\}$ and
$\mathbb A_+=\{t=x+\varepsilon y\in\mathbb A\ | \ t\bar t>0, x>0\}$.
The para-complex projective space $\mathbb AP^m$ is defined
 to be the quotient of $\mathbb A_0^{m+1}$ under the
equivalence relation: $Z\simeq W$ if $Z=tW$ for some $t\in\mathbb
A_+$ (see \cite{erdem,gadea-amilibia}).

  We note that  $H^1= \{t\in \mathbb A_+\ |\
t\bar t=1\}$.
 The restriction of the projection $\mathbb A_0^{m+1}\to
\mathbb AP^m=\mathbb A_0^{m+1}/\mathbb A_+$  to $H^{2m+1}_{m}$,
gives the Hopf submersion
\begin{eqnarray}
\pi_{\mathbb A}::H^{2m+1}_{m}\to\mathbb AP^m=H^{2m+1}_{m}/H^1.
\end{eqnarray}
Moreover, there exists a unique pseudo-Riemannian metric $g'$ on
$\mathbb AP^m$ such that $\pi_{\mathbb A}$
 is a pseudo-Riemannian
submersion with totally geodesic fibres \cite{erdem}.
The space $(\mathbb AP^m,g')$ is a complete para-holomorphic space form and its
para-holomorphic curvature is $-4$.

The construction of $\mathbb BP^m$ is analogous to the para-complex
projective space. We have
\begin{eqnarray}
\mathbb BP^m= \{z\in\mathbb B^{m+1}\ |\ \langle
z,z\rangle=1\}/\{t\in \mathbb B|\ t\bar t=1
 \}=H^{4m+3}_{2m+1}/H^3_1,
\end{eqnarray} and there exists a unique pseudo-Riemannian metric
$g'$ on
 $\mathbb BP^m$ such that the projection
\begin{eqnarray}
\pi_{\mathbb B}:H^{4m+3}_{2m+1}\to
 \mathbb BP^m=H^{4m+3}_{2m+1}/H^3_1
\end{eqnarray}
 is a
pseudo-Riemannian submersion with totally geodesic fibres
\cite{blazic}. Moreover, $(\mathbb BP^m,g')$ is a complete, simply
connected, para-quaternionic space form of para-quaternionic
curvature $-4$ (see \cite{blazic}).

\subsection{The Hopf pseudo-Riemannian submersions between pseudo-hyperbolic spaces: the Hopf construction}
All Hopf pseudo-Riemannian submersions between (real)
pseudo-hyperbolic spaces can explicitly be obtained by the Hopf
construction.

A bilinear map $G:\mathbb R^p\times\mathbb R^q\to\mathbb R^n$ is
said to be an orthogonal multiplication if $G$ is norm-preserving,
that is $|G(x,y)|=|x| |y|$ for any $x\in\mathbb R^p, y\in\mathbb R^q$
(see \cite{baird-wood,smith}). A {\it Hopf construction} is a map
$\varphi:\mathbb R^p\times\mathbb R^q\to\mathbb R^{n+1}$ defined by
$\varphi(x,y)=(|x|^2-|y|^2,2G(x,y))$, for some orthogonal
multiplication  $G$ (see \cite{baird-wood,smith}). The Hopf
construction can provide several examples of harmonic morphisms (see
\cite{konderak,smith}), and we would like to refer the reader to the
 beautiful book \cite{baird-wood} due to Baird and Wood for other nice results
on this topic.
 Since the sectional curvatures $K$, $K'$ of the total and base spaces of any pseudo-Riemannian
submersion between real space forms must obey  $K'=4K$, we are
forced to consider the map $\varphi(x,y)/2$ instead.

Let $F\in\{\mathbb C, \mathbb A, \mathbb H, \mathbb B, \mathbb O,
\mathbb O'\},$  and  let $\varphi_1,\varphi_2:F\times F\to \mathbb
R\times F$ be, respectively, the maps given by
\begin{eqnarray}
 \varphi_1(x,y)=((|x|^2-|y|^2)/2,\bar xy) \ \ \ \text{and}\ \
 \varphi_2(x,y)=((|x|^2+|y|^2)/2,\bar xy)
\end{eqnarray}
for any $x,y\in F$, where $\bar x$ denotes the conjugation of $x$ in
$F$ and as usual $|x|^2=x\bar x$, $|y|^2=y\bar y$. For convenience,
we denote $t_1=(|x|^2-|y|^2)/2\in\mathbb R$,
$t_2=(|x|^2+|y|^2)/2\in\mathbb R$ and $w=\bar xy\in F$. Since
 $|w|^2=|\bar xy|^2=|x|^2|y|^2$ for any $x,y\in F$, it is easy to see that
\begin{itemize}
\item[(i)] if $|x|^2+|y|^2=1$, then $t_1^2+|w|^2=1/4$;
\item[(ii)] if $|x|^2-|y|^2=1$, then $t_2^2-|w|^2=1/4$.
\end{itemize}
Setting  $d=\dim_{\mathbb R}F$, we identify $F\times F\simeq\mathbb
R^{2d}$ via
\begin{eqnarray}\label{e:ident}
((x^1,\cdots,x^d),(y^1,\cdots,y^d))\simeq(x^1,y^1,\cdots,x^d,y^d).
\end{eqnarray}

When $F\in\{\mathbb C,  \mathbb H,   \mathbb O\}$, we consider the
following restrictions of $\varphi_1$ and  $\varphi_2$ to
$H^{2d-1}_{2d-1}$ and to $H^{2d-1}_{d-1}$, respectively:
 $$\varphi_1:H^{2d-1}_{2d-1}=\{(x,y)\in F^2 \ |\
|x|^2+|y|^2=1\}\to H^{d}_{d}(-4)=\{(t_1,w)\in\mathbb R\times F \ |\
 t_1^2+|w|^2=1/4\},$$
$$\varphi_2:H^{2d-1}_{d-1}=\{(x,y)\in F^2 \ |\ |x|^2-|y|^2=1\}\to H^{d}(-4)
=\{(t_2,w)\in\mathbb R\times F \ |\ t_2^2-|w|^2=1/4\}.$$ This simple
construction gives six Hopf pseudo-Riemannian submersions with
totally geodesic fibres:
$$\pi_1:H^3_3\to H^2_2(-4)=\mathbb CH^1_1,\ \ \
  \pi_2:H^7_7\to H^4_4(-4)=\mathbb HH^1_1,\ \ \
  \pi_3:H^{15}_{15}\to H^8_8(-4),
$$
$$
 \pi_4:H^3_1\to H^2(-4)=\mathbb CH^1,\ \ \
 \pi_5:H^7_3\to H^4(-4)=\mathbb HH^1, \ \ \
 \pi_6:H^{15}_{7}\to H^8(-4)=\mathbb OH^1.
$$
The first three submersions are the well-known Hopf fibrations
between spheres.

When $F\in\{\mathbb A,  \mathbb B, \mathbb O'\}$, the restriction of
$\varphi_1$ to $H^{2d-1}_{d-1}$,
$$\varphi_1:H^{2d-1}_{d-1}=\{(x,y)\in F^2 \ |\ |x|^2+|y|^2=1\}\to H^{d}_{d/2}(-4)
=\{(t_1,w)\in\mathbb R\times F \ |\  t_1^2+|w|^2=1/4\},$$ gives
another three Hopf pseudo-Riemannian submersions with totally
geodesic fibres between pseudo-hyperbolic spaces:
$$
  \pi_7:H^3_1\to H^2_1(-4)=\mathbb AH^1, \ \ \
  \pi_8:H^7_3\to H^4_2(-4)=\mathbb BH^1, \ \ \
  \pi_9:H^{15}_{7}\to H^8_4(-4).
$$
Note that, for $F\in\{\mathbb A,  \mathbb B, \mathbb O'\}$, the
restriction of $\varphi_2$ to $H^{2d-1}_{d-1}$ will give the same
$\pi_7,\pi_8,\pi_9$.
 In \cite{konderak}, Konderak  constructed the harmonic
 morphisms $2\pi_7$ and $2\pi_8$ via the Hopf construction (see also
\cite[Examples 14.6.5-6]{baird-wood}). For identification
\eqref{e:ident} of $\mathbb O'\times\mathbb O' \simeq\mathbb
R^{16}$, the Hopf pseudo-Riemannian submersion $\pi_9:H^{15}_{7}\to
H^8_4(-4)$  can be  written explicitly as
\begin{eqnarray*}
\pi_9(x_1,y_1,\cdots,x_8,y_8)&=&\big( (x_1^2 + x_2^2 + x_3^2 + x_4^2
- x_5^2 -
 x_6^2 - x_7^2 - x_8^2
      -y_1^2 - y_2^2- y_3^2 - y_4^2+ y_5^2\\&&
  \hspace*{-3cm}    + y_6^2 + y_7^2 + y_8^2)/2,
     \ \  x_1y_1 + x_2  y_2 + x_3  y_3 + x_4  y_4 - x_5  y_5 - x_6  y_6 - x_7  y_7 -
        x_8 y_8, \, \, -x_2 y_1
\\&&  \hspace*{-3cm}
          + x_1  y_2 + x_4 y_3 - x_3  y_4 - x_6  y_5 + x_5 y_6 +
        x_8 y_7 - x_7 y_8, \ \
         -x_3y_1 - x_4y_2 + x_1  y_3 + x_2 y_4
\\&& \hspace*{-3cm}
          - x_7 y_5 -
        x_8 y_6 + x_5 y_7 + x_6  y_8,
\ \
         -x_4 y_1 + x_3  y_2 - x_2 y_3 + x_1 y_4 -
        x_8 y_5 + x_7 y_6 - x_6  y_7
\\&& \hspace*{-3cm}
        + x_5  y_8,
\ \
   -x_5 y_1 - x_6 y_2 - x_7 y_3 - x_8 y_4 + x_1y_5 + x_2y_6 +
        x_3 y_7 + x_4 y_8, \ \
        -x_6y_1 + x_5 y_2
\\&& \hspace*{-3cm}
          - x_8 y_3 + x_7 y_4 - x_2 y_5 +
        x_1 y_6 - x_4 y_7 + x_3 y_8,
\ \
         -x_7  y_1 + x_8 y_2 + x_5 y_3 - x_6 y_4 -
        x_3 y_5
\\&& \hspace*{-3cm}
        + x_4 y_6 + x_1 y_7 - x_2 y_8,
\ \
         -x_8 y_1 - x_7  y_2 + x_6 y_3 +
        x_5 y_4 - x_4 y_5 - x_3 y_6 + x_2 y_7 + x_1 y_8\big).
\end{eqnarray*}

Note that $\pi_1$, $\pi_2$, $\pi_4$, $\pi_5$, $\pi_7$, $\pi_8$
fall, respectively, in the  categories $\pi_{\mathbb C}$, $\pi_{\mathbb H}$ ,
$\pi_{\mathbb C}$, $\pi_{\mathbb H}$, $\pi_{\mathbb A}$,
$\pi_{\mathbb B}$ of \S\ref{s:3.1} and \S\ref{s:3.2}. Define
$\pi^1_{\mathbb O}=\pi_3$, $\pi^2_{\mathbb O}=\pi_6$, $\pi_{\mathbb
O'}=\pi_9$. To the best of our knowledge $\pi_{\mathbb O'}$  does
not appear in the literature.

  The construction of the Hopf
pseudo-Riemannian submersions solves the existence problem for each
class $(a)$-$(g)$ of Theorem \ref{t:main}. In the following sections, we
approach the uniqueness.

\begin{remark}
 The Hopf pseudo-Riemannian submersions are homogeneous, i.e. of the
 form  $\pi:G/K\to G/H$ with
 $K\subset H$ closed Lie subgroups:\\
 $\pi_{\mathbb C}:H^{2m+1}_{2t+1}=SU(m-t,t+1)/SU(m-t,t)\to
       \mathbb CH^m_t=SU(m-t,t+1)/S(U(1)U(m-t,t)),$
\\
 $\pi_{\mathbb H}:H^{4m+3}_{4t+3}=Sp(m-t,t+1)/Sp(m-t,t)\to
     \mathbb HH^m_t=Sp(m-t,t+1)/Sp(1)Sp(m-t,t),$
\\
    $\pi_{\mathbb A}:H^{2m+1}_{m}=SU(m+1, \mathbb A)/SU(m,\mathbb A)\to\mathbb
    AP^m=SU(m+1,\mathbb A)/S(U(1,\mathbb A)U(m,\mathbb A))$,\
\\
    $\pi_{\mathbb B}:H^{4m+3}_{2m+1}=Sp(m+1,\mathbb B)/Sp(m,\mathbb B)\to
    \mathbb BP^m=Sp(m+1,\mathbb B)/Sp(1,\mathbb B)Sp(m,\mathbb B),$
\\
    $\pi^1_{\mathbb O}:H^{15}_{15}=Spin(9)/Spin(7)\to
     H^8_8(-4)=Spin(9)/Spin(8),$
\\
     $\pi^2_{\mathbb O}:H^{15}_7=Spin(8,1)/Spin(7)\to
     H^8(-4)=Spin(8,1)/Spin(8),$
\\
    $\pi_{\mathbb O'}:H^{15}_{7}=(Spin(5,4)/Spin(3,4))_0\to H^8_{4}(-4)=(Spin(5,4)/Spin(4,4))_0$.
\\
 By Harvey's book \cite[p.~312]{harvey-book}, each of
$Spin(5,4)/Spin(3,4)$ and $Spin(5,4)/Spin(4,4)$ has two connected
components: a pseudo-sphere and a pseudo-hyperbolic space. Here
$(\cdot)_0$ denotes the pseudo-hyperbolic component.

By analogy to Hopf Riemannian submersions from spheres \cite{bes},
each of the canonical variations of   $\pi_{\mathbb B}$,
$\pi_{\mathbb H}$, $\pi^1_{\mathbb O}$, $\pi^2_{\mathbb O}$ and
$\pi_{\mathbb O'}$ gives a new homogeneous Einstein metric on the
pseudo-hyperbolic space. The classification problem of homogeneous
Einstein metrics on pseudo-hyperbolic spaces shall be discussed
somewhere else.
\end{remark}

\section{The geometry of the base space}\label{s4}

An important step of the proof of Theorem \ref{t:main} is to
establish that the base space is either a real space form or a
special Osserman pseudo-Riemannian manifold. By the classification
of complete, simply connected, special Osserman pseudo-Riemannian
manifolds
\cite{garcia-rio-kupeli-vazquez-lorenzo,special-osserman-classification},
we explicitly get the geometry of the base space, and then we see
that the dimensions and the indices  of the total space and of the
base are those claimed in Theorem \ref{t:main}. First, we recall
Proposition 3.8 from \cite{tohoku}, which provides the completeness
and the simply-connectedness of the base space.
 \bpr\label{p:5}
    Let $\pi:M\to B$ be a pseudo-Riemannian submersion
    with connected totally geodesic fibres from a complete connected
    pseudo-Riemannian manifold $M$ onto a pseudo-Riemannian manifold $B$. Then
    $B$ is complete. Moreover, if $M$ is simply connected, then $B$
    is also simply connected.
\epr
 Let $\pi:M\to B$ be a pseudo-Riemannian submersion.
 We  use the following notation throughout the paper:
 $n=\dim B$, $s=\mathrm{index}\, B$, $F_{b}=\pi^{-1}(b)$ for some $b\in B$,
 $r=\dim F_{b}$ and $r'=\mathrm{index}\, F_{b}$.
\subsection{The construction of a special basis $\mathcal B$ of $\mathcal
H$ along a fibre} A key ingredient for understanding the geometry of
the base and of the fibres is the construction of a special
orthonormal basis $\mathcal B$ of $\mathcal H$ along a fibre, which
we recall from \cite{tohoku}. First, we state the following lemma,
which provides useful properties of
  O'Neill's integrability tensor for a constant curvature total
space.
 \ble[\cite{tohoku}]\label{lem:2}
    Let $\pi:M\to B$ be a pseudo-Riemannian submersion
    with connected totally geodesic
    fibres from a pseudo-Riemannian manifold $M$ with
    constant curvature
    $c\neq 0$. Then the following assertions are true$:$
\begin{itemize}
    \item[(a)]
    If $X$ is a horizontal vector such that $g(X,X)\neq 0$, then
    the map $A_X:\mathcal V\to\mathcal H$
        given by $A_X(V)=A_XV$
        is injective and the map
    $A^*_X:\mathcal{H}\to\mathcal{V}$ given by
    $A^*_X(Y)=A_XY$ is surjective.
\item[(b)] If $X$, $Y$ are the horizontal lifts
    along the fibre
    $\pi^{-1}(\pi(p))$, $p\in M$, of
    two vectors $X', Y'\in T_{\pi(p)}B$, respectively, $g'(X',X')\neq 0$
    and $(A_XY)(p)=0$, then $A_XY=0$ along the fibre
    $\pi^{-1}(\pi(p))$.
\end{itemize}
\ele The proof of Lemma \ref{lem:2} relies on   the O'Neill
equations. Corollary \ref{p:oneill2}(b) simply gives
\begin{eqnarray}\label{e:AXAXV}
  A^*_XA_XV=-c g(X,X)V
\end{eqnarray}
 for every
vertical vector field $V$, which implies (i). By Corollary
\ref{p:oneill2}(a), we get (ii).

 Let
$p\in M$ and let
      $\{v_{1p},\dots  ,v_{rp}\}$ be an orthonormal basis in $\mathcal
      V_p$.
      Let $X'\in T_{\pi(p)}B$ such that $g'(X',X')=\pm 1$
      and let $X$ be the horizontal lift along the fibre
      $\pi^{-1}(\pi(p))$  of $X'$. Let $Y_1,Y_2,\dots,Y_r$ be the horizontal lifts along the fibre
      $\pi^{-1}(\pi(p))$ of
      $$\frac{1}{cg(X,X)}\pi_*A_Xv_{1p}, \ \frac{1}{cg(X,X)}\pi_*A_Xv_{2p} ,\dots,
      \frac{1}{cg(X,X)}\pi_*A_Xv_{rp},$$, respectively.
      For each $i\in\{1,\dots,r\}$, we consider the vector
      $v_i=A_XY_i$ defined along the fibre $\pi^{-1}(\pi(p))$.
By Corollary \ref{p:oneill2}(a),
       $\{v_1,v_2,\dots,v_r\}$ is an orthonormal basis of $\mathcal
       V_q$
      at any $q\in\pi^{-1}(\pi(p))$  (see \cite{tohoku}), which can
      be restated as the following lemma.

\ble[\cite{tohoku}]\label{l:o1}
     In the set-up of Lemma \ref{lem:2}, the fibres are
     parallelizable.
\ele

      Set $L_0=X$.
      For every  integer $\al$ with $1\leq\al <n/(r+1)$,
      let $L_\al$ be a horizontal vector field along the fibre
      $\pi^{-1}(\pi(p))$ such that
\begin{itemize}
\item[(1)]      $L_\al$ is
      the horizontal lift of some unit vector
      (i.e., $g(L_\al,L_\al)\in\{-1,1\}$), and
\item[(2)]
      $L_\al$ is orthogonal to
      $L_0,L_1,\dots  ,L_{\al-1}$, and
\begin{eqnarray}\label{e:LAAA}
      L_\al(p)\in\ker A^*_{L_0(p)}\cap\ker A^*_{L_1(p)}\cap\dots
      \cap\ker A^*_{L_{\al-1}(p)}.
\end{eqnarray}
\end{itemize}
   Condition \eqref{e:LAAA} is nothing but the statement that  $L_\al(p)$
   is orthogonal to any vector in the system $\{L_0(p),A_{L_0}v_1(p),\cdots,A_{L_0}v_r(p),\cdots,
             L_{\al-1}(p),A_{L_{\al-1}}v_1(p),\cdots,A_{L_{\al-1}}v_r(p)\}$.
      Moreover, by Lemma \ref{lem:2}(b), $L_\al(q)$ belongs to
      $\ker A^*_{L_0(q)}\cap\ker A^*_{L_1(q)}\cap\dots
      \cap\ker A^*_{L_{\al-1}(q)}$ for every $q\in\pi^{-1}(\pi(p))$.
     In the set-up of
     Lemma \ref{lem:2},
     Proposition
     \ref{p:oneill}(c) implies that
\begin{eqnarray}
\mathcal B=\{L_0,A_{L_0}v_1,\cdots,A_{L_0}v_r,\cdots,
             L_{k-1},A_{L_{k-1}}v_1,\cdots,A_{L_{k-1}}v_r\}
\end{eqnarray}
is an orthonormal basis of $\mathcal H_q$ for any
$q\in\pi^{-1}(\pi(p))$ (see \cite{tohoku}). It is worth pointing out
that any element in $\mathcal B$ is basic along the fibre
$\pi^{-1}(\pi(p))$ by \eqref{e:LAAA} and Proposition
\ref{p:oneill}(a) (see \cite{tohoku}). Such a basis $\mathcal B$ is
said to be a {\it special basis}.

Counting   the time-like vectors of $\mathcal B$, we get the
following proposition.

\bpr[\cite{tohoku}]\label{p:3}
    In the set-up of Lemma \ref{lem:2}, we have
    $n=k(r+1)$ for some  positive integer $k$
    and $s=q_1(r'+1)+q_2(r-r')$ for
    some  non-negative integers $q_1$, $q_2$ with
    $q_1+q_2=k$.
\end{proposition}

The following corollary will be needed later.
\begin{corollary}[\cite{tohoku}]\label{c:son}
If $s\in\{0,n\}$, then $r'=r$ (i.e. the metrics induced on fibres
are negative definite).
\end{corollary}

We now split the problem of identifying the geometry of $B$ into two
cases: $(i)$ $n=r+1$ (that is, $k=1$), and $(ii)$   $n\not=r+1$ (that is,
$k>1$).

\subsection{Case $n=r+1$}
This case features a constant curvature base space:
\begin{proposition}\label{t:const}
 In the set-up of Theorem \ref{t:main},
  $n=r+1$ if and only if $B$ has constant curvature $-4$.
\end{proposition}
\begin{proof}
Let $b\in B$, $X'\in T_bB$ such that $g'(X',X')=\pm 1$ and $p\in
\pi^{-1}(b)$. Let $X\in\mathcal H_p$ be the horizontal lift of $X'$.

Assuming $n=r+1$, that is,
 $\dim\mathcal H_p=\dim\mathcal V_p+1$, we see that $A_X:\mathcal V_p\to X^\perp=\{Y\in\mathcal H_p\ | \ g(X,Y)=0\}$
 is bijective, and thus, for every $Y\in X^\perp$, we can write $Y=A_XV$ for
 some vertical vector $V$. By \eqref{e:AXAXV}, we get
\begin{eqnarray}\label{e:o1}
 g(A_XY, A_XY)=g(A_XA_XV,A_XA_XV)=g(X,X)^2g(V,V).
\end{eqnarray}
 On the other hand, by Corollary \ref{p:oneill2}(b), we have
\begin{eqnarray}\label{e:o2}
 g(Y,Y)=g(A_XV,A_XV)=-g(X,X)g(V,V).
\end{eqnarray}
 Combining equations
 \eqref{e:o1} and \eqref{e:o2}, we simply get
 $g(A_XY,A_XY)=-g(X,X)g(Y,Y)$ for every $Y\in X^\perp$, which
 implies that $A_XA_XZ=g(X,X)Z-g(X,Z)X$ for any horizontal vector
 $Z$. Now, by Corollary \ref{p:oneill2}(a), we  obtain
\begin{eqnarray}\label{e:o3}
R'(X,Y,X,Y)&=&-g(X,X)g(Y,Y)+g(X,Y)^2+3g(A_XY,A_XY)\nonumber\\
 &=&-4\big(g(X,X)g(Y,Y)-g(X,Y)^2\big),
\end{eqnarray}
which means that  $B$  has constant curvature $-4$.

 Conversely, if $B$ has constant curvature $-4$, then,  by
\eqref{e:o3}, we get  $g(A_XY,A_XY)=\\-g(X,X)g(Y,Y)$ for every $Y\in
X^\perp$, which implies $A_XA_XY=g(X,X)Y$ for every $Y\in X^\perp$.
 Therefore, by \eqref{e:AXAXV}, $A_X:\mathcal V\to X^\perp$ is bijective
 with its inverse is given by $(A_X)^{-1}(Y)=(1/(g(X,X)))A_XY$, for $Y\in
 X^\perp$.
 As a consequence,
 $n-1=\dim X^\perp=\dim \mathcal V_p=r$.
\end{proof}

\begin{theorem}\label{t:constant0}
In the set-up of Theorem \ref{t:main}, if  $n=r+1$ and $0<s<n$,
  then $\pi$ falls into one of the following cases:
\begin{itemize}
\item[(a)] $\pi:H^3_1\to H^2_1(-4)=\mathbb AH^1$,
\item[(b)] $\pi:H^7_3\to H^4_2(-4)=\mathbb BH^1$,
\item[(c)] $\pi:H^{15}_7\to H^8_4(-4)$.
\end{itemize}
\end{theorem}

\begin{proof}
 First, recall that  $B$ has
constant curvature $-4$ by Proposition \ref{t:const}. Let $X$,
$Y\in\mathcal H_p$ such that $g(X,X)=1$ and $g(Y,Y)=-1$.
 Let
$\mathcal B=\{X,A_Xv_1,\cdots, A_Xv_r\}$, $\mathcal
B'=\{Y,A_Yv'_1,\cdots, A_Yv'_r\}$ be two special bases of $\mathcal
H_p$.  The index of $\mathcal B$, the number of time-like vectors,
 is $r-r'$, while the index of $\mathcal B'$ is
$r'+1$. Therefore, $r=2r'+1$, $s=r'+1$, and $n=2(r'+1)$. The
pseudo-Riemannian submersion $\pi$ is of the form
  $\pi:H^{4r'+3}_{2r'+1}\to B^{2r'+2}_{r'+1}$.

  By a theorem due to Reckziegel \cite{rec},
  the horizontal distribution $\mathcal H$  of a pseudo-Riemannian submersion with totally geodesic fibres
  is an Ehresmann connection, and thus, by Ehresmann \cite{ehr},  $\pi$ is a locally trivial fibration,
  which always comes with a long exact homotopy sequence
\begin{eqnarray}\label{e:lehs}
\cdots\to\pi_2(B)\to\pi_1(F_{\pi(p)})\to\pi_1(H^{4r'+3}_{2r'+1})\to\pi_1(B)\to\pi_0(F_{\pi(p)})\to\cdots.
\end{eqnarray}

Now, we proceed in two cases: (i) $r'=0$ and (ii) $r'\geq 1$.

 {\bf Case  $r'=0$}.  Since the fibres are connected, totally
 geodesic, one-dimensional submanifolds (when $r'=0$),
  any fibre is the image of a space-like geodesic in
 $H^{4r'+3}_{2r'+1}$. Thus, the fibres are diffeomorphic to the real
 line (see \cite[p.~113]{onei}) and $\pi_1(F_{\pi(p)})=0$. The long exact homotopy
 sequence \eqref{e:lehs} gives $\pi_1(B)=\pi_1(H^3_1)=\mathbb Z$.
  Because  $B$ is of constant curvature $-4$,  and, by
 Proposition
 \ref{p:5}, is also complete, it simply follows that $B$ is isometric to the
 pseudo-hyperbolic space $H^2_1(-4)$, and that corresponds to (a).

  {\bf Case $r'\geq 1$}.  By the long exact homotopy
 sequence \eqref{e:lehs}, and by
 $\pi_1(H^{4r'+3}_{2r'+1})=\pi_1(S^{2r'+1})=0$, we get $\pi_1(B)=0$.
 The manifold  $B$
 is additionally complete  and of constant curvature $-4$. Therefore
 $B$ must be isometric to $H^{2r'+2}_{r'+1}(-4)$. The case $r'=1$
 corresponds to (b).

  We now assume that $r'\geq 2$.
  Since, for $r'\geq 2$,
  $\pi_2(B)=\pi_2(H^{2r'+2}_{r'+1}(-4))=\pi_2(S^{r'+1}\times\mathbb R^{r'+1})=0$   and
  $\pi_1(H^{4r'+3}_{2r'+1})=\pi_1(S^{2r'+1}\times\mathbb R^{2r'+2})=0$,
  the long exact homotopy sequence \eqref{e:lehs} gives
  $\pi_1(F_{\pi(p)})=0$.
  On the other hand, since the fibres are totally geodesic in $H^{4r'+3}_{2r'+1}$,  the
  fibres are complete and of curvature $-1$. Therefore, the fibres
  must be isometric to $H^{2r'+1}_{r'}$. By Lemma \ref{l:o1}, the fibres
  are also parallelizable, and that restricts the choices of $r'\geq 2$ to $r'\in\{3,7\}$. The value
  $r'=3$ corresponds to the cases (c).

  We now show that the case $r'=7$ is not possible, namely we see that
  there is no
   pseudo-Riemannian submersion $\pi:H^{31}_{15}\to H^{16}_8(-4)$ with
  connected totally geodesic fibres. By Ranjan
  \cite{ran}, the linear map
  $\mathcal U: \mathcal V_p\to\mathrm{Hom}(\mathcal H_p,\mathcal H_p)$ given by $\mathcal
  U(V)(X)=A_XV$ extends to a Clifford representation
  $\mathcal U: Cl(\mathcal V_p,-\hat g)\to\mathrm{Hom}(\mathcal H_p,\mathcal
  H_p)$, namely
  $\mathcal U(v )\mathcal U(w)+\mathcal U(w)\mathcal
  U(v)=2g(v,w)\mathrm{Id}$ for every $v,w\in \mathcal V_p$, because of
  Corollary \ref{p:oneill2}(b).
   This makes the sixteen-dimensional space $\mathcal H_p$ a $Cl(\mathcal V_p)$-module,
  which, as usual, decomposes into
  irreducible $Cl(\mathcal V_p)$-modules. On the other hand,
  the signature of the inner product $-\hat g(v,w)=-g(v,w)$
  of $\mathcal V_p$ is  $(7,8)$, and from the
  Classification Table of the Clifford algebras \cite[p.~29]{michelsohn},
  we see that
  $Cl(\mathcal V_p, -\hat g)=Cl_{(7,8)}=\mathbb R(128)\oplus\mathbb R(128)$.
  In consequence,  any irreducible $Cl(\mathcal V_p)$-module
  is of dimension 128, and thus the dimension of $\mathcal H_p$ is too small
  to allow a nontrivial Clifford representation
  $\mathcal U: Cl(\mathcal V_p)\to\mathrm{Hom}(H_p,H_p)$ as above.
\end{proof}

The case $s=0$ corresponds to a Riemannian base space which was
completely classified in \cite{dga}, while the case $s=n$ is of a
Riemannian submersion from spheres (classified in \cite{esc,ran})
when we apply a change of signs of the metrics of the total and of
the base spaces. By Corollary \ref{c:son}, the metrics induced on
fibres are negative definite if $s\in\{0,n\}$.
\begin{theorem}[\cite{dga,esc,ran}]\label{t:constant3}
In the set-up of Theorem \ref{t:main},  we assume $n=r+1$.  Then the
following assertions are true$:$
\begin{itemize}
\item[$(i)$] If $s=0$,
  then $\pi$ is one of the following:
\begin{itemize}
 \item[$(a)$] $\pi:H^3_1\to H^2(-4)$,
 \ \ \ \ $(b)$  $\pi:H^7_3\to H^4 (-4)$,
 \ \ \ \ $(c)$
$\pi:H^{15}_7\to H^8(-4)$.
\end{itemize}
\item[$(ii)$]  If $s=n$,
  then $\pi$ is one of the following:
\begin{itemize}
\item[$(a')$] $\pi:H^3_3\to H^2_2(-4)$,
  \ \ \ \ $(b')$ $\pi:H^7_7\to H^4_4 (-4)$,
  \ \ \ \ $(c')$ $\pi:H^{15}_{15}\to H^8_8(-4)$.
\end{itemize}
\end{itemize}
\end{theorem}

\subsection{Case $n\not= r+1$}

We show that $B$ is a complete, simply connected, special Osserman
pseudo-Riemannian manifold.

\subsubsection{Special Osserman manifolds} Following
\cite{garcia-rio-kupeli-vazquez-lorenzo}, we recall the definitions
of  a Jacobi operator and of a special Osserman pseudo-Riemannian
manifold.
\begin{defn}
Let $(B,g')$ be a pseudo-Riemannian manifold and let $R'$ be the
Riemann curvature tensor of $(B,g')$. For $x\in T_bB$, we consider
the linear map $R'(\cdot,x)x:T_bB\to T_bB.$ Since
$g'(R'(z,x)x,x)=0$, we have $\mathrm{Im}( R'(\cdot,x)x)\subset
x^{\perp}$, where $x^\perp=\{y\in T_bB\ | \ g'(y,x)=0\}$. For $x\in
S_bB=\{x\in T_bB\ |\ g'(x,x)=\pm 1\}$, the restriction
$R'_x:x^\perp\to x^\perp$ of $R'(\cdot,x)x$ to $x^\perp$ is called
the {\it Jacobi operator} with respect to $x$, that is,
$R'_x(z)=R'(z,x)x$ for $z\in x^\perp$.
\end{defn}
\begin{defn}\label{d:definition-special-osserman}
A pseudo-Riemannian manifold $(B,g')$ is called {\it special
Osserman} if the following two conditions are satisfied at each
$b\in B$:
\begin{itemize}
\item[(I)] For every $x\in S_bB$ the Jacobi operator $R'_x:x^\perp\to x^\perp$ is
diagonalizable with exactly two distinct eigenvalues
$\varepsilon_x\lambda$ and $\varepsilon_x\mu$, where
$\varepsilon_x=g'(x,x)$ and $\lambda,\mu\in\mathbb R$.
 \item[(II)] Let
 $E_\lambda(x)=\mathrm{span}\{x\}\oplus\ker(R'_x-\varepsilon_x\lambda
 Id)$. For each $x\in S_bB$, if $z\in E_\lambda(x)\cap S_bB$,
  then $E_\lambda(x)=E_\lambda(z)$, and moreover if
 $y\in S_bB\cap\ker(R'_x-\varepsilon_x\mu Id)$, then
 $x\in \ker(R'_y-\varepsilon_y\mu Id)$.
\end{itemize}
\end{defn}
The values $\lambda$ and $\mu$ involved in the previous definition
are not interchangeable, for example if  $(B,g',J)$ is the complex
or the para-complex pseudo-hyperbolic space of real dimension
$2n>2$, then $\mu=\lambda/4$ and $\ker(R'_x-\varepsilon_x\lambda
Id)=\mathrm{span}\{Jx\}$ is one-dimensional, while
$\ker(R'_x-\varepsilon_x\mu Id)=\{x,Jx\}^\perp=\{z\ |\
g'(z,x)=g'(z,Jx)=0\}$ is $(2n-2)$-dimensional.

\subsubsection{The base space is Special Osserman} For a
pseudo-Riemannian submersion $\pi:(M,g)\to (B,g')$, we denote by
$R'_{X'}$ the Jacobi operator of $(B,g')$ with respect to a vector
$X'\in T_bB$ and for   $X,Y\in\mathcal H_p$ we also denote by
$R'_{X}Y$ the horizontal lift of $R'_{\pi_{*}(X)}(\pi_*Y)$ and we
consider $R'_{X}$ as an operator $R'_{X}:X^\perp\to X^\perp$, with
$X^\perp=\{Y\in\mathcal H_p \ |\ g(Y,X)=0\}$.

\begin{theorem}\label{t:osserman}
In the set-up of Theorem \ref{t:main}, if $n\not=r+1$, then $B$ is
special Osserman.
\end{theorem}
\begin{proof}
Let $b\in B$, $X'\in S_bB$, $Z'\in T_bB$ and $p\in \pi^{-1}(b)$. Let
$X, Z\in\mathcal H_p$ be the horizontal lifts of $X'$ and $Z'$, respectively. By Corollary \ref{p:oneill2}(a),  $R'_{X}$ is given by
\begin{eqnarray}\label{e:Rprim}
R'_{X}(Z)=R'(Z,X)X=R(Z,X)X-3A_XA_XZ=R_XZ-3A_XA_XZ.
\end{eqnarray}
Let $\{v_1,v_2,\cdots, v_r\}$ be an orthonormal basis in $\mathcal
V_p$, that is, $g(v_i,v_j)=\varepsilon_i\delta_{i,j}$ with
$\varepsilon_i\in\{-1,1\}$. Let
$$\mathcal B=\{L_0,A_{L_0}v_1,\cdots, A_{L_0}v_r, \cdots, L_{k-1},A_{L_{k-1}}v_1,\cdots,
A_{L_{k-1}}v_r\}$$ be a special basis of $\mathcal H_p$, that is an
orthonormal basis $\mathcal B $ with $L_0=X$ and
$A_{L_\alpha}L_\beta=0$ for every $\alpha,\beta\in\{0,\cdots,k-1\}$.
We show that $R'_{X}:X^\perp\to X^\perp$ is diagonalizable with
respect to $\mathcal B$ and $R'_{X}$ has exactly two eigenvalues. By
\eqref{e:Rprim} and \eqref{e:AXAXV}, we have
 \begin{eqnarray}\label{e:o4}
 \ \ \ \ \ \
 R'_{X}(A_Xv_i)&=&R_X(A_Xv_i)-3A_XA_XA_Xv_i\nonumber\\
  &=&-g(X,X)A_Xv_i-3g(X,X)A_Xv_i
   = -4\varepsilon_XA_Xv_i,
 \end{eqnarray}
which gives
$R'_{X'}(\pi_*(A_Xv_i))=\pi_*(R'_{X}(A_Xv_i))=-4\varepsilon_{X'}\pi_*(A_Xv_i)$.
 Since  $$0=g(A_Xv_j,A_{L_\alpha}
v_i)=-g(v_j,A_XA_{L_\alpha} v_i)$$ for every $i,j$ and every
$\alpha\geq 1$, we get
 $A_XA_{L_\alpha} v_i=0$, which implies that
 \begin{eqnarray}\label{e:o5}
 R'_{X }( A_{L_\alpha} v_i)&=&R_X(A_{L_\alpha}v_i)-3A_XA_XA_{L_\alpha}v_i
  =-g(X,X)A_{L_\alpha}v_i
   = -\varepsilon_XA_{L_\alpha}v_i .
 \end{eqnarray}
Projecting \eqref{e:o5} to the base, we have
$R'_{X'}(\pi_*(A_{L_\alpha} v_i))=
-\varepsilon_{X'}\pi_*(A_{L_\alpha}v_i)$.
  Since
$A_XL_\alpha=0$ by construction, we see that
\begin{eqnarray}\label{e:o6}
 \ \ \  \ \ \  \ \ \  \ \ \  R'_{X}( L_\alpha) =R_X(L_\alpha)-3A_XA_X L_\alpha
  =-g(X,X)L_\alpha=-\varepsilon_XL_\alpha
\end{eqnarray}
for every $\alpha\geq 1$ and every $i$. Therefore
$R'_{X'}(\pi_*(L_\alpha)) =-\varepsilon_{X'}\pi_*(L_\alpha).$
 Summarizing, the
Jacobi operator $R'_{X'}$ is diagonalizable with the eigenvalues
$-4\varepsilon_{X'}$ and $-\varepsilon_{X'}$, and moreover their
eigenspaces are:
\begin{eqnarray}\label{e:R'X'}
 \ker(R'_{X'}+4\varepsilon_{X'}\mathrm{Id})&=&\{\pi_*(A_Xv_1),\cdots,\pi_*(A_Xv_r)\}\text{
and, }\\
\ker(R'_{X'}+\varepsilon_{X'}\mathrm{Id})&=&\{\pi_*(L_1),\pi_*(A_{L_1}v_1),\cdots,\pi_*(A_{L_1}v_r),\cdots,\nonumber\\
&&
\pi_*(L_{k-1}),\pi_*(A_{L_{k-1}}v_1),\cdots,\pi_*(A_{L_{k-1}}v_r)\}.
 \end{eqnarray}
Now, we check that Condition (II) of Definition
\ref{d:definition-special-osserman} holds.

\begin{lemma}\label{l:1}
If $Y'\in E_{-4}(X')$, $g'(X',X')=\pm 1$ and $g'(Y',Y')=\pm 1$, then
$X'\in E_{-4}(Y')$.
\end{lemma}

\begin{proof}[Proof of Lemma \ref{l:1}]
By \eqref{e:R'X'},
 $$E_{-4}(X')=\mathrm{span}\{X'\}\oplus \ker(R'_{X'}+4\varepsilon_{X'}\mathrm{Id})
  =\mathrm{span}\{\pi_*X,\pi_*(A_Xv_1),\cdots,\pi_*(A_Xv_r)\},$$ and, thus, the
  horizontal lift $Y$ of $Y'$ satisfies
\begin{eqnarray}\label{e:Y-X}
  Y=aX+A_XU
\end{eqnarray}
   for some
  $a\in\mathbb R$ and some vertical vector $U$. By \eqref{e:Y-X},
\begin{eqnarray}\label{e:new}
 g(A_XU,A_XU)=g(Y,Y)-a^2g(X,X).
\end{eqnarray}
  To prove $X'\in
  E_{-4}(Y')$, it is sufficient to show that  $X$ can be written as
\begin{eqnarray}\label{e:X-Y}
  X=bY+A_YW
 \end{eqnarray}
   for some $b\in\mathbb R$ and some vertical vector $W$. Applying
   $A_Y$ to \eqref{e:X-Y}, we get $A_YX=bA_YY+A_YA_YW=g(Y,Y)W$,
   which gives  $W=-A_XY/(g(Y,Y))$. Similarly, applying $A_X$ to
   \eqref{e:Y-X}, we obtain $A_XY=A_XA_XU=g(X,X)U$.
   Substituting $Y$ and $W$ into
   \eqref{e:X-Y}, we obtain an equation in  $b\in\mathbb R$
\begin{eqnarray}\label{e:EQ}
   X&=&b(aX+A_XU)-\frac{g(X,X)}{g(Y,Y)}A_{aX+A_XU}U \text{,  which is equivalent to}\\
   X&=&baX-\frac{g(X,X)}{g(Y,Y)}A_{A_XU}U+(b-\frac{ag(X,X)}{g(Y,Y)})A_XU.\label{e:EQ22}
 \end{eqnarray}
By
Corollary \ref{p:oneill2}(b),
\begin{eqnarray}\label{e:op}
g(A_XU,A_ZU)=-g(X,Z)g(U,U)
 \end{eqnarray}
 for every horizontal
vectors $X,Z$ and for every vertical vector $U$. Since $A$ is
 skew-symmetric with respect to $g$ and alternating, we have
  $g(A_XU,A_ZU)=-g(A_ZA_XU,U)=g(A_{A_XU}Z,U)=-g(Z,A_{A_XU}U)$, which
  by \eqref{e:op},
  implies that $A_{A_XU}U=g(U,U)X$. Then
\begin{eqnarray*}
baX-\frac{g(X,X)}{g(Y,Y)}A_{A_XU}U&=&(ba-\frac{g(X,X)g(U,U)}{g(Y,Y)})X=(ba+\frac{g(A_XU,A_XU)}{g(Y,Y)})X\\
&=&(ba+\frac{g(Y,Y)-a^2g(X,X)}{g(Y,Y)})X=X-a(b-\frac{ag(X,X)}{g(Y,Y)})X,
 \end{eqnarray*}
by \eqref{e:new}. Therefore, \eqref{e:EQ22} has the unique solution
$b=\frac{ag(X,X)}{g(Y,Y)}$.
\end{proof}

\begin{lemma}\label{l:2}
If $Y'\in \ker(R'_{X'}+\varepsilon_{X'}\mathrm{Id})$, $g'(X',X')=\pm
1$  and $g'(Y',Y')=\pm 1$, then
$X'\in\ker(R'_{Y'}+\varepsilon_{Y'}\mathrm{Id})$.
\end{lemma}

\begin{proof}[Proof of Lemma \ref{l:2}] Let $X$ and $Y$ be the
horizontal lifts of $X'$ and $Y'$, respectively. The Jacobi operator $R'_{X'}$
satisfies
\begin{eqnarray}
R'_{X'}(Y')=\pi_*(R_X(Y)-3A_XA_XY)=-g'(X',X')Y'-3\pi_*(A_XA_XY)
\end{eqnarray}
 for any $Y'\in
X'^\perp$. Therefore, $Y'\in
\ker(R'_{X'}+\varepsilon_{X'}\mathrm{Id})$ if and only if
$A_XA_XY=0$. Since, by Lemma \ref{lem:2}(a), $A_X:\mathcal
V\to\mathcal H$ is injective,   $A_XY=0$, hence, $A_YX=0$, which
implies  that
$R'_{Y'}(X')=\pi_*(-3A_YA_YX+R_Y(X))=-g'(Y',Y')X'=-\varepsilon_{Y'}X'$.
\end{proof}
These conclude that $B$ is a special Osserman pseudo-Riemannian
manifold.
\end{proof}

In the next theorem, we identify the geometry of the base space and
we find the dimension and the index of the total space in terms of
the geometry of the base space.
\begin{theorem}\label{t:constant2}
 Let $\pi:H^{n+r}_{s+r'}\to B^{n}_{s}$ be a pseudo-Riemannian
under the assumptions of Theorem \ref{t:main}. If  $n\not=r+1$
  then $\pi$ falls in one of the following cases:
\begin{itemize}
    \item[$(a)$]  $H^{2m+1}_{2t+1}\to\mathbb CH^m_t$,\ \
    \ \ \    \ \ \ $(b)$   $H^{2m+1}_{m}\to\mathbb AH^m$,
    \item[$(c)$]  $H^{4m+3}_{4t+3}\to\mathbb HH^m_t$,\ \
   \ \ \     \ \ \ $(d)$   $H^{4m+3}_{2m+1}\to\mathbb BH^m$,
    \item[$(e)$]  $H^{23}_{7}\to\mathbb OH^2$,
    \ \ \    \ \ \   $(f)$   $H^{23}_{15}\to\mathbb OH^2_1$,
     \ \ \    \ \ \   $(g)$    $H^{23}_{23}\to\mathbb OH^2_2$,
     \ \ \    \ \ \  $(h)$   $H^{23}_{q}\to\mathbb O'P^2$,
\end{itemize}
for $0\leq t\leq m$ and $m\geq 2$, and for some $8\leq q\leq 15$.
\end{theorem}
\begin{proof}
We first prove that $B$ is simply connected. When ${s+r'}>1$,
$H^{n+r}_{s+r'}$ is simply connected and thus, by Proposition
\ref{p:5}, $B$ is also simply connected. If $s+r'=1$, then either
$(i)$ $s=0$ and $r'=1$, or $(ii)$ $s=1$ and $r'=0$.

In the case $(i)$ $s=0$ and $r'=1$, the base space is Riemannian,
which,  by Magid \cite{mag}, must be isometric to $\mathbb
 CH^m$, and thus $B$ is simply connected.

In the case $(ii)$ $s=1$ and $r'=0$, $B$ is Lorentzian Osserman at
the point $p$, which by Garc\'\i a-R\'\i o, Kupeli and V\' azquez-Lorenzo \cite{garcia-rio-kupeli-vazquez-lorenzo}, it
must be of constant curvature at the point $p$. On the other
hand, $B$ has constant curvature if and only if $n=r+1$. This
contradicts our working assumption $n\not=r+1$. These conclude that
$B$ is simply connected.

By the classification theorem of simply connected, complete  special
Osserman pseudo-Riemannian manifolds
\cite{special-osserman-classification,garcia-rio-kupeli-vazquez-lorenzo},
$B$ is isometric to one of the following:
\begin{itemize}
\item[(a)]  a  definite or indefinite complex space form of signature
$(2m-2s,2s)$, $0\leq s\leq m$;
\item[(b)]  a  definite or  indefinite quaternionic space form of signature
$(4m-4s,4s)$, $0\leq s\leq m$;
\item[(c)] a para-complex space form of signature $(m, m)$;
\item[(d)] a para-quaternionic space form  of signature
$(2m, 2m)$;
\item[(e)] a Cayley plane of  octonions with definite or indefinite
metric, or a Cayley plane of para-octonions with indefinite metric
of signature $(8,8)$.
\end{itemize}
Any non-flat complete, simply connected, para-complex space form is
isometric to the symmetric space $SL(m+1,\mathbb R)/(SL(m,\mathbb
R)\times\mathbb R)=\mathbb AP^m$ (see
\cite{special-osserman-classification,cruceanu,garcia-rio-kupeli-vazquez-lorenzo}),
 and  any non-flat complete, simply connected para-quaternionic
space form is isometric to the symmetric space $Sp(m+1,\mathbb
B)/(Sp(m,\mathbb B)Sp(1,\mathbb B))=Sp(2m+2,\mathbb
R)/(Sp(2m,\mathbb R)SL(2,\mathbb R))=\mathbb BP^m$ (see
\cite{special-osserman-classification,dancer-jorgensen-swann,garcia-rio-kupeli-vazquez-lorenzo,GMV}).

By the proof of Theorem \ref{t:osserman}, the values $\lambda$ and
$\mu$ of Definition \ref{d:definition-special-osserman} are
negative, namely $\lambda=-4$ and $\mu=-1$. Then $B$ must be
isometric to one of the following spaces:
\begin{eqnarray}\label{e:basespace} \ \ \ \  \ \ \ \
\mathbb CH^m_t,\ \mathbb HH^m_t,\ \mathbb AP^m, \ \mathbb BP^m,\
\mathbb OH^2,\ \mathbb OH^2_1,\ \mathbb OH^2_2,\  \text{or }\
\mathbb O'P^2,
\end{eqnarray}
with $m\geq 2$ and $0\leq t\leq m$. By \eqref{e:R'X'}, we simply
have  $\dim\ker(R'_{X'}+4\varepsilon_{X'}\mathrm{Id})=r$ , and in
particular the following conditions are satisfied.
\begin{itemize}
\item[(a)] If $B\in\{\mathbb CH^m_t, \mathbb AP^m\}$, then
$\ker(R'_{X'}+4\varepsilon_{X'}\mathrm{Id})=\mathrm{span}\{IX'\}$,
where $I$ is a complex or para-complex structure. Thus $r=1$ and
$n+r=2m+1$.
\item[(b)] If $B\in\{\mathbb HH^m_t, \mathbb BP^m\}$, then
$\ker(R'_{X'}+4\varepsilon_{X'}\mathrm{Id})=\mathrm{span}\{IX',JX',KX'\}$,
with $\{I,J,K\}$ a local quaternionic or para-quaternionic
structure. Therefore, $r= 3$ and $n+r=4m+3$.
\item[(c)] If $B\in\{\mathbb OH^2_i, \mathbb O'P^2\}_{0\leq i\leq 2}$, then
$\dim\ker(R'_{X'}+4\varepsilon_{X'}\mathrm{Id})=7$. Thus $r=7$ and
$n+r=23$.
\end{itemize}

Now, we   find the index of the total space for the choices of $B$
in \eqref{e:basespace}.

{\bf Case 1: $B \in\{\mathbb CH^m_t, \mathbb HH^m_t, \mathbb
OH^2_i\}_{ 0\leq t\leq m,\, 0\leq i\leq 2}$.} In this case, the
Riemann tensor satisfies
\begin{eqnarray}\label{e:R'}
R'(X',Y',X',Y')\leq -(g(X',X')g(Y',Y')-g(X',Y')^2)
\end{eqnarray}
 for any $X'$, $Y'$ vectors on $B$. Let $\{v_i\}_{i\in\{1,\cdots,r\}}$ be an
orthonormal basis of $\mathcal V_p$ and let $X$ be the horizontal lift
of a non-null vector $X'\in T_{\pi(p)}B$. Taking $Y'=\pi_*(A_Xv_i)$,
inequality \eqref{e:R'} becomes
\begin{eqnarray}\label{e:R'2}
R'(\pi_*X ,\pi_*(A_Xv_i),\pi_*X ,\pi_*(A_Xv_i))\leq -g(X
,X)g(A_Xv_i,A_Xv_i).
\end{eqnarray}
On the other hand   by Corollary \ref{p:oneill2}(a) and by
\eqref{e:AXAXV},
$$R'(\pi_*X ,\pi_*(A_Xv_i),\pi_*X ,\pi_*(A_Xv_i))=
 -4g(X,X)g(A_Xv_i,A_Xv_i).$$ Now, \eqref{e:R'2} implies $0\leq
 3g(X,X)g(A_Xv_i,A_Xv_i)=-g(X,X)^2g(v_i,v_i)$ for any $i$. Thus, the fibres are
 negative definite. Therefore, in Case 1, $\pi$ should be in one
 of  (a), (c), (e)-(g) of Theorem \ref{t:constant2}.
     Note that, in Case 1, $B$ is isotropic which means that
      for any
      $b\in B$ and any $t\in\mathbb R$, the group of isometries
      of $B$ preserving $b$ acts transitively on $\{Z\in T_bB\, |\,  g'(Z,Z)=t,\ Z\not=0 \}$
       (see \cite[p.~367]{wolf}).

 {\bf Case 2: $B=\mathbb AP^m$.} Since $B=\mathbb AP^m$ is a
para-quaternionic space form of para-holomorphic curvature
$\lambda=-4$,
\begin{eqnarray}\label{e:R'para}
R'(X',Y',X',Y')\geq -(g(X',X')g(Y',Y')-g(X',Y')^2).
\end{eqnarray}
By a similar argument to Case 1, specializing \eqref{e:R'para} for a
non-null vector $X'$ and $\pi_*(A_Xv_1)$ we get $0\geq
 3g(X,X)g(A_Xv_1,A_Xv_1)=-g(X,X)^2g(v_1,v_1)$ and thus the fibres are
 positive definite and $\pi$ falls in (b).

{\bf Case 3: $B=\mathbb BP^m$.} We shall show that the fibres have
signature $(2,1)$. Note that $(\mathbb BP^m,g')$ has a natural
para-quaternionic K\"ahler structure and its curvature tensor
satisfies the relation
\begin{eqnarray}\label{e:R'3}
 R'(X',Y',X',Y')&=&-(g'(X',X')g'(Y',Y')-g'(X',Y')^2 \nonumber\\
 &&-3g'(J_1X',Y')^2-3g'(J_2X',Y')^2+3g'(J_3X',Y')^2),
\end{eqnarray}
where $\{J_1,J_2,J_3\}$ is a local para-quaternionic structure, a
triple of $(1,1)$-tensors
 satisfying
 $J_1J_2=-J_2J_1=J_3$,
$J_i^2=\varepsilon_i\mathrm{Id}$, $g'(J_iX',Y')+g'(X',J_iY')=0$ and
$\varepsilon_1=\varepsilon_2=-\varepsilon_3=1$. Obviously, for any
$X',Y'$ such that $g'(J_3X',Y')=0$ we have
\begin{eqnarray}\label{e:R'3inequality}
 R'(X',Y',X',Y')\geq -(g'(X',X')g'(Y',Y')-g'(X',Y')^2).
\end{eqnarray}
Let $X'\in T_b\mathbb BP^m$ such that $g'(X',X')=\pm 1$ and let $X$
be its horizontal lift at $p\in\pi^{-1}(b)$. Let $J_3X\in\mathcal
H_p$ be the horizontal lift of $J_3X'$. By \eqref{e:R'3},
$$R'(X',J_3X',X',J_3X')=-4g'(X',X')g'(J_3X',J_3X'),$$ and thus
$$g(A_XJ_3X,A_XJ_3X)=-g(X,X)g(J_3X,J_3X)=-g(X,X)^2=-1,$$
by Corollary \ref{p:oneill2}(a). Let $\{v_1,v_2,v_3\}$ be an
orthonormal basis of $\mathcal V_p$ such that $v_3=A_XJ_3X$. We
simply note that $g(v_3,v_3)=-1$. For $i\in\{1,2\}$, taking
$Y'=\pi_*(A_Xv_i)$ in \eqref{e:R'3inequality}, we get
\begin{eqnarray}\label{e:R'5}
 R'(X',\pi_*(A_Xv_i),X',\pi_*(A_Xv_i))\geq
 -g'(X',X')g'(\pi_*(A_Xv_i),\pi_*(A_Xv_i)).
\end{eqnarray}
On the other hand, $R'(X,A_Xv_i,X,A_Xv_i)=-4g(X,X)g(A_Xv_i,A_Xv_i)$.
Thus, \eqref{e:R'5} becomes $0\geq
3g(X,X)g(A_Xv_i,A_Xv_i)=-3g(X,X)^2g(v_i,v_i)$ for $i\in\{1,2\}$.
Therefore, $g(v_i,v_i)>0$ for $i\in\{1,2\}$.
 \end{proof}

To see that the cases (e)-(h) of Theorem \ref{t:constant2} never
occur, we first recall the notion of Clifford structure.
\subsubsection{Clifford structures}

We adapt the definition of Clifford structure introduced by Gilkey
\cite{gilkey} and Gilkey, Swann and Vanhecke
\cite{gilkey-swann-vanhecke} to pseudo-Riemannian geometry.

\begin{defn}\label{d:def-cliff}
Let $(B,g')$ be a pseudo-Riemannian manifold and let $R'$ be its
curvature tensor.
 The space $(B,g')$ has a $\mathrm{Cliff}(\nu)$-structure if
at every point $b$ there exist (1,1)-tensors $J_1,J_2,\cdots, J_\nu
$ such that
\begin{eqnarray}\label{e:def-cliff}
  \ \ \ \  \ \ R'(x,y)z&=&
 \lambda_0(g'(y,z)x-g'(x,z)y)+\frac{1}{3}\sum_{s=1}^{\nu}\varepsilon_s(\lambda_s-\lambda_0)
  (g'( J_sy,z)J_sx
  \nonumber\\
 && -g'(J_sx,z)J_sy
  - 2g'(J_sx,y)J_sz
  ),
\end{eqnarray}
for any $x,y,z\in T_bB$, where
$\lambda_0,\lambda_1,\cdots,\lambda_\nu:B\to\mathbb R$,
$\lambda_s(b)\not=\lambda_0(b)$ for $s\geq 1$, and
 $g'(J_sx,y)=-g'(x, J_sy)$ and
 $J_sJ_t+J_tJ_s=-2\varepsilon_s\delta_{s,t}\mathrm{Id}$, with $\varepsilon_s=\pm
 1$.
\end{defn}

The Jacobi operator at the point $b$ of a manifold with a
$\mathrm{Cliff}(\nu)$-structure   is given by:
\begin{eqnarray}\label{e:jac-cliff}
 R'_y(x)= \lambda_0g'(y,y) x
+ \sum_{s=1}^{\nu}\varepsilon_s(\lambda_s-\lambda_0)
  g'(x,J_sy)J_sy,
\end{eqnarray}
for any $x\in y^\perp$. Moreover,
\begin{eqnarray}\label{e:clifford}
R'_y(J_sy)&=&\lambda_sg'(y,y)J_sy \text{ for any }
 s\in\{1,\cdots,\nu\}\text{ and}\\
R'_y(x)&=& \lambda_0g'(y,y)x \text{
 for any } x\in\{y,J_1y,\cdots,J_ry\}^\perp,\label{e:clifford2}
\end{eqnarray}
 and thus  a pseudo-Riemannian manifold with
  a $\mathrm{Cliff}(\nu)$-structure
 is  pointwise Osserman (see \cite{gilkey-book}).

In the Riemannian setup, Clifford structures turned out to be a
very valuable tool for the Osserman Conjecture. In
\cite{gilkey-swann-vanhecke}, Gilkey, Swann and Vanhecke  suggested
a two-step approach: (i) show that the pointwise Osserman condition
implies the existence of a Clifford structure with
\eqref{e:clifford}, \eqref{e:clifford2}, and (ii) find the manifolds
having the curvature tensors of (i). Using this approach,
Nikolayevsky proved the Osserman conjecture in dimension $n\not=16$;
see
 \cite{nik2,nik}. In dimension n=16,
 the Cayley planes $\mathbb OH^2$, $\mathbb
OP^2$ do not admit Clifford structures %\cite{nik3},
\cite[p.~510]{nik} and the Osserman Conjecture remains open.

 Since the curvature tensor formulae of the Cayley planes of
octonions or of para-octonions are similar to that of $\mathbb
OP^2$, in particular
 the eigenspace of the Jacobi
operator for $\lambda=-4$ satisfies
\begin{eqnarray}\label{e:eigen}
 \ \  \ \ \  \ker(R'_{(a,b)}+4\varepsilon_{(a,b)}\mathrm{Id})=\left\{
\begin{array}{cl}
\big\{\big(c,\ \frac{1}{|a|^2}(b\bar a)c\big)\ |\
\mathrm{Re}\,(c\bar a)=0\big\},
&\text{ if } |a|^2\not=0,\\
\big\{\big(\frac{1}{|b|^2}(a\bar b)d,\ d\big)\ |\
\mathrm{Re}\,(d\bar b)=0\big\}, &\text{ if } |b|^2\not=0,
 \end{array}
\right.
\end{eqnarray}
for any $(a,b)\in S_bB$ (see \cite{held-stavrov-vankoten}), one can
easily see,  by analogy to \cite[p.~510]{nik}, that $ \mathbb
OH^2_2, \mathbb OH^2_1, \mathbb OH^2$, $\mathbb O'P^2$ do not admit
Cliff(7)-structures. To exclude (e)-(h) of Theorem
\ref{t:constant2}, it is now sufficient to establish the following
theorem.

\begin{theorem}\label{t:cliffordstructure}
Let $\pi:M\to B$ be a pseudo-Riemannian submersion with connected
totally geodesic fibres. If $M$ has constant curvature $c\not=0$,
then $B$ has a $\mathrm{Cliff}(r)$-structure.
\end{theorem}

\begin{proof}
Without loss of the generality,  we may assume $c=\pm 1$. Let $p\in
M$ and $b=\pi(p)\in B$. Let $\{v_1,\cdots, v_r\}$ be an orthonormal
basis of $\mathcal V_p$.
 For any $1\leq s\leq r$, let
$\varepsilon_s=cg(v_s,v_s)\in\{-1,1\}$ and let $J_s(X')=\pi_*(A_X
v_s)$ where $X\in T_pM$ is the horizontal lift of $X'\in T_bB$. For
any vertical  vector $v\in\mathcal V_p$, we define the linear map
$A^v:\mathcal H_p\to\mathcal H_p$ given by
 $A^v(x)=A_xv$ for $x\in\mathcal H_p$.
 Since $M$ has constant curvature $c$, by Ranjan's
 paper \cite{ran},
 we have
\begin{eqnarray}\label{e:ran-clif}
A^vA^w+A^wA^v=-2cg(v,w)Id,
\end{eqnarray}
 for any $v,w$ vertical vectors. Thus
 $J_sJ_t+J_tJ_s=-2cg(v_s,v_t)\mathrm{Id}=-2\varepsilon_s\delta_{s,t}\mathrm{Id}$.
 Also, by Ranjan's paper \cite{ran}, we have $g(A^vX,Y)=-g(X,A^vY)$ for any $X,Y\in\mathcal H_p$,
 which simply implies $g'(J_sX',Y')=-g'(X',J_sY')$ for every $X',Y'\in T_bB$ and every $1\leq s\leq r$.

Now, we show that the Jacobi operator of $B$ satisfies
\eqref{e:jac-cliff}. Let $X', Y'\in T_bB$ with $g'(Y',Y')=\pm 1$, and
$g(X',Y')=0$. Let $X$ and $Y$ be the horizontal lifts of $X'$ and $Y'$, respectively.
Let
$$\mathcal B=\{L_0,A_{L_0}v_1,\cdots, A_{L_0}v_r, \cdots, L_{k-1},A_{L_{k-1}}v_1,\cdots,
A_{L_{k-1}}v_r\}$$ be a special basis of $\mathcal H_p$  such that
$L_0=Y$. We recall that $\mathcal B$ is orthonormal and that
$A_{L_\alpha}L_\beta=0$ for every $\alpha,\beta\in\{0,\cdots,k-1\}$,
by construction.
 $X$ can be written as
\begin{eqnarray}\label{e:Xwrites}
X&=&g(X,Y)Y
 +\sum_\alpha\frac{g(X,L_\alpha)}{g(L_\alpha,L_\alpha)}L_\alpha\nonumber\\
 &&+\sum_i\frac{g(X,A_Yv_i)}{cg(Y,Y)g(v_i,v_i)}A_Yv_i
 +\sum_{i,\alpha}\frac{g(X,A_{L_\alpha}
 v_i)}{cg(L_\alpha,L_\alpha)g(v_i,v_i)}A_{L_\alpha}v_i.
\end{eqnarray}
Since $\mathcal B$ is orthonormal, $A_YA_{L_\alpha}v_i=0$ by the proof
of Theorem \ref{t:osserman}. Applying $A_YA_Y$ to \eqref{e:Xwrites},
we get
 $$
 A_YA_YX=\sum_i\frac{g(X,A_Yv_i)}{cg(Y,Y)g(v_i,v_i)}A_YA_YA_Yv_i=-c\sum_i \varepsilon_ig(X,A_Yv_i)
 A_Yv_i=-\sum_i \varepsilon_icg(X,J_iY)J_iY
 $$
Then 
\begin{eqnarray}\label{e:jacobiop}
\ \ \ \ \ \ \ R'_{Y'}(X')
=\pi_*(R_YX-3A_YA_YX)=cg'(Y',Y')X'+3c\sum_i
\varepsilon_ig'(X',J_iY')J_iY'.
\end{eqnarray}
Polarizing  \eqref{e:jacobiop}, we get
\begin{eqnarray*}
  R'(X',Y')Z'&=&
 c(g'(Y',Z')X'-g'( X',Z') Y')\\
 && +c\sum_{i=1}^{r}\varepsilon_i
  (g'(J_iY',Z') J_iX'
  -g'(J_iX',Z') J_iY'
  - 2g'(J_iX',Y') J_iZ'
  ).
\end{eqnarray*}
\end{proof}

\begin{corollary}\label{t:nocayley}
There  are no pseudo-Riemannian submersions $\pi:H^{23}_t\to B$ with
connected totally geodesic fibres from a $23$-dimensional
pseudo-hyperbolic space $H^{23}_t$ onto any of the Cayley
pseudo-hyperbolic planes of octonions  $ \mathbb OH^2_2, \mathbb
OH^2_1, \mathbb OH^2$, or onto the Cayley projective plane of
para-octonions $\mathbb O'P^2$.
\end{corollary}

\begin{remark}
Ranjan \cite{ran} proved that there are no Riemannian submersions
$\pi:S^{23}\to\mathbb OP^2$ with connected, totally geodesic fibres
(that is, (g) of Theorem \ref{t:constant2}). For a topological proof
of this fact we refer the reader to \cite{tang}.
\end{remark}

\section{The Theorem of Uniqueness}\label{s5}
To prove Theorem \ref{t:main} we need the following Theorem of
Uniqueness.
\begin{theorem}\label{t:criterium}
Let $\pi_1,\pi_2:H^{a}_{l}\to B$ be two pseudo-Riemannian
submersions with connected totally geodesic fibres from a
pseudo-hyperbolic space onto a pseudo-Riemannian manifold. Then
there exists an isometry $f:H^{a}_{l}\to H^{a}_{l}$ such that
$\pi_2\circ f=\pi_1$. In particular, $\pi_1$ and $\pi_2$ are
equivalent.
\end{theorem}
\begin{proof}
The main ideas of the proof are: (1) for a given basepoint $b$
construct special bases $\mathcal B^1$ and $\mathcal B^2$ for the
fibres $F^1_b$ and $F^2_b$, respectively, such that $\mathcal B^1$
and $\mathcal B^2$ have the same projections to the base B
and (2) show that the unique
isometry sending $\mathcal B^1$ into $\mathcal B^2$ preserves the
integrability tensors everywhere and sends fibres into fibres.

Let $b\in B$ and $p,q\in H^{a}_{l}$ such that $\pi_1(p)=\pi_2(q)=b$.
%For simplicity, we set $b=\pi_1(p)=\pi_2(q)$.
We denote by $\mathcal
V^1$ and $\mathcal V^2$ the vertical distributions of $\pi_1$  and
$\pi_2$, and by $\mathcal H^1$ and  $\mathcal H^2$  the horizontal
distributions of $\pi_1$ and $\pi_2$, respectively.

 Let $\{v_{1p},\cdots,v_{rp}\}$ be an orthonormal basis of $\mathcal V^1_p$ and let $X'\in
 T_bB$ such that $g'(X',X')=\pm 1$.
 We  denote by $X^1$ and $X^2$  the $\pi_1$- and  $\pi_2$-horizontal lifts of $X'$
 along the fibres $F^1_{b}=\pi_1^{-1}(b)$  and
   $F^2_{b}=\pi_2^{-1}(b)$, respectively.
 Let $(Y^1_1,Y^1_2,\dots,Y^1_r)$ and $(Y^2_1,Y^2_2,\dots,Y^2_r)$ be the $\pi_1$- and $\pi_2$-horizontal lifts
       of
      $$\big(\frac{1}{-g(X^1,X^1)}\pi_{1*}A_{X^1}v_{1p}, \ \frac{1}{-g({X^1},{X^1})}\pi_{1*}A_{X^1}v_{2p} ,\dots,
      \frac{1}{-g({X^1},{X^1})}\pi_{1*}A_{X^1}v_{rp}\big),$$
      along $F^1_{b}$ and $F^2_{b}$, respectively.
      For each $i\in\{1,\dots,r\}$, we consider the vectors
      $v_i=A^1_{X^1}Y^1_i$,   defined along  $F^1_{b}$,
       and  $w_i=A^2_{X^2}Y^2_i$ along  $F^2_{b}$.
     By \S\ref{s4}, $\{v_1,\cdots,v_r\}$ is a global orthonormal
     basis of vector fields on $F^1_b$, and we claim that so is
   $\{w_1,\cdots,w_r\}$. Indeed, by Corollary \ref{p:oneill2}(a), we
   see that
\begin{eqnarray*}
g(w_i,w_j)&=&g(A^2_{X^2}Y^2_i,A^2_{X^2}Y^2_j)\\
&=&(1/3)(R'(\pi_{2*}X^2,\pi_{2*}Y^2_j,\pi_{2*}X^2,\pi_{2*}Y^2_j)
 -g(X^2, X^2)g(Y^2, Y^2)+g(X^2,Y^2)^2)\\
&=&(1/3)(R'(\pi_{1*}X^1,\pi_{1*}Y^1_j,\pi_{1*}X^1,\pi_{1*}Y^1_j)
 -g(X^1, X^1)g(Y^1, Y^1)+g(X^1,Y^1)^2)\\
&=&g(A^1_{X^1}Y^1_i,A^1_{X^1}Y^1_j)=g(v_i,v_j)=\varepsilon_i\delta_{ij}
\end{eqnarray*}
along $F^2_{b}$. Let $\mathcal
B^1=\{L^1_0,A^1_{L^1_0}v_1,\cdots,A^1_{L^1_0}v_r,\cdots,
      L^1_{k-1},A^1_{L^1_{k-1}}v_1,\cdots,A^1_{L^1_{k-1}}v_r\}$
be a special basis of  $\mathcal H^1$   along $F^1_b$
  such that $L^1_0=X^1$ (and $A^1_{L^1_\alpha}L^1_\beta=0$). Let $L^2_1,\cdots, L^2_{k-1}$ be the
  $\pi_2$-horizontal lifts of $\pi_{1*}L^1_1,\cdots,
  \pi_{1*}L^1_{k-1}$
along $F^2_b$. We take  $L^2_0=X^2$. Let    $$\mathcal
B^2=\{L^2_0,A^2_{L^2_0}w_1,\cdots,A^2_{L^2_0}w_r,\cdots,
      L^2_{k-1}, A^2_{L^2_{k-1}}w_1,\cdots,A^2_{L^2_{k-1}}w_r\}.$$

\begin{lemma}\label{l:basicA1A2}
(i) The vector field $A^2_{X^2}w_i$ is basic along $F^2_b$ and
$\pi_{1*}A^1_{X^1}v_i=\pi_{2*}A^2_{X^2}w_i$,  for every $i$.

(ii) We have $A^2_{X^2}L^2_\alpha=0$ and $A^2_{L^2_\alpha}L^2_{\beta}=0$ for
every $\alpha$ and $\beta$.

(iii) The basis $\mathcal B^2$ is a special basis of $\mathcal H^2$ along
$F^2_b$ and
$\pi_{1*}A^1_{L^1_{\alpha}}v_i=\pi_{2*}A^2_{L^2_{\alpha}}w_i$, for
every $i$ and $\alpha$.
\end{lemma}
\begin{proof} Let $Z'\in T_bB$, and let $Z^1$ and $Z^2$ be the $\pi_1$- and
$\pi_2$-horizontal lifts of $Z'$ along $F^1_b$ and $F^2_b$,
respectively. By Corollary \ref{p:oneill2}(a), we get
\begin{eqnarray*}
g(A^2_{X^2}w_i, Z^2)&=&
  -g(A^2_{X^2}Y^2_i,A^2_{X^2}
  Z^2)=(1/3)(R(X^2,Y^2_i,X^2,Z^2)-R'(X^2,Y^2_i,X^2,Z^2))\\
  &=&
  (1/3)(R(X^1,Y^1_i,X^1,Z^1)-R'(X^1,Y^1_i,X^1,Z^1))
  =-g(A^1_{X^1}Y^1_i,A^1_{X^1}Z^2)\\
&=& g(A^1_{X^1}v_i, Z^1),
\end{eqnarray*}
which simply  implies (i). By (i), we see that
\begin{eqnarray}\label{e:A222}
g(A^2_{X^2}L^2_\alpha, w_i)&=&
  -g(L^2_\alpha, A^2_{X^2}w_i)=-g'(\pi_{2*}L^2_\alpha, \pi_{2*}A^2_{X^2}w_i)
  \nonumber\\ &=&-g'(\pi_{1*}L^1_\alpha, \pi_{1*}A^1_{X^1}v_i)=g(A^1_{X^1}L^1_\alpha,
  v_i)=0,
\end{eqnarray}
 for every $i$ and $\alpha$. Thus, $A^2_{X^2}L^2_\alpha=0$.
Therefore, by Proposition \ref{p:oneill}(a), we obtain that
\begin{eqnarray*}
2g(A^2_{L^2_\alpha}L^2_{\beta},w_i)&=&2g(A^2_{L^2_\alpha}L^2_{\beta},A^2_{X^2}Y^2_i)
=R'( L^2_\alpha,L^2_{\beta},X^2,Y^2_i)\\
&&
 -R(L^2_\alpha,L^2_{\beta},X^2,Y^2_i)
 +g(A^2_{L^2_{\beta}}{X^2},A^2_{L^2_\alpha}{Y^2_i})
 -g(A^2_{L^2_\alpha}{X^2},A^2_{L^2_{\beta}}{Y^2_i})\\
&& =R'(L^1_\alpha,L^1_{\beta},X^1,Y^1_i)
 -R(L^1_\alpha,L^1_{\beta},X^1,Y^1_i)=2g(A^1_{L^1_\alpha}L^1_{\beta},v_i)=0,
\end{eqnarray*}
for every $i$. Thus $A^2_{L^2_\alpha}L^2_{\beta}=0$ and hence
$\mathcal B^2$ is a special basis of $\mathcal H^2$.

  By Proposition \ref{p:oneill}(c), $A^2_{L^2_{\alpha}}w_i$ is
basic along $F^2_b$ (for details see \cite[Lemma 3.4]{tohoku}), and
by an argument similar to \cite[Lemma 3.4]{tohoku} one can see that
$\pi_{1*}A^1_{L^1_{\alpha}}v_i=\pi_{2*}A^2_{L^2_{\alpha}}w_i$.
\end{proof}

Since $\mathcal B^1$ and $\mathcal B^2$ are special bases, they are
orthonormal, by \S\ref{s4}.
 Let $F:T_pH^a_l\to
T_qH^a_l$ be the linear isometry given by $F(v_i)=w_i$,
$F(A^1_{L^1_\alpha}v_i)=A^2_{L^2_\alpha}w_i$,
$F(L^1_\alpha)=L^2_\alpha$, for any $1\leq i\leq r$, $0\leq
\alpha\leq k-1$.  Since $H^{a}_{l}$
 is a frame-homogeneous space, there exists an isometry $f:H^a_l\to
 H^a_l $ such that $f(p)=q$ and $f_{*p}=F$   (see \cite{onei,wolf}).
It remains to prove that $\pi_2\circ f=\pi_1$.

We say that the condition $(\star)$ is satisfied at $x\in H^a_l$ if
\begin{eqnarray*}
(\star) \ \ \pi_2(f(x))=\pi_1(x), \ \ \ f_{*x}(\mathcal
H^1_x)=\mathcal H^2_{f(x)}, \ \ \ f_*(A^1_EF)=A^2_{f_*E}f_*F \text{
for any } E,F\in T_xH^a_l.
\end{eqnarray*}

We will proceed in four steps.
\begin{itemize}
\item[{\bf Step 1.}] $(\star)$ holds at $p$.
\item[{\bf Step 2.}] $(\star)$ holds at every  $z\in F^1_{b}$.
\item[{\bf Step 3.}] If $\tilde\gamma:[0,1]\to H^a_l$ is a
               $\pi_1-$horizontal geodesic  with $\tilde\gamma(0)\in F^1_{b}$,
                then $(\star)$ holds at any point
  $\tilde\gamma(t)$, where $t\in[0,1]$.
\item[{\bf Step 4.}]  $\pi_2(f(x))=\pi_1(x)$ for any $x\in H^a_l$.
\end{itemize}

{\bf Proof of Step 1.} From the definition of $F$, we simply have
$\pi_2(f(p))=\pi_1(p)$ and
\begin{eqnarray}\label{e:step1}
f_{*p}(\mathcal H^1_p)=\mathcal H^2_{f(p)}.
\end{eqnarray}
 We recall that the
vectors  of $\mathcal B^1$
 are basic along $F^1_{b}$. Since
\begin{eqnarray}
A^1_{A^1_{L^1_\al}v_i}A^1_{L^1_\bet}
v_j=g(L^1_\al,L^1_\bet)\hat\na^1_{v_i}v_j
\end{eqnarray}
along $F^1_{b}$ (see \cite{tohoku}) and since $A^1$ is alternating,
we see that $\hat\na^1_{v_i}v_j=(1/2)[v_i,v_j]$. Similar relations
hold for
  $\pi_2$, and, at  $p$, we simply have
$f_*[v_i,v_j]=[f_*v_i,f_*v_j]=[w_i,w_j]$. Therefore,
\begin{eqnarray}\label{e:A1A2}
f_*(A^1_{A^1_{L^1_\al}v_i}A^1_{L^1_\bet}
v_j)=A^2_{f_*(A^1_{L^1_\al}v_i)}f_*(A^1_{L^1_\bet} v_j).
\end{eqnarray}
By the definition of $f$ and \eqref{e:A1A2}, we get
$f_{*p}(A^1_EF)=A^2_{f_{*p}E}f_{*p}F$ for any $E,F\in T_pH^a_l$.

{\bf Proof of Step 2.} The following lemma  shall be needed right
away.

\begin{lemma}[{\cite[p.~105]{onei}}]\label{l:esc-oneill}
Let $N_1, N_2$ be two complete,  connected, totally geodesic
pseudo-Riemannian submanifolds of a pseudo-Riemannian manifold $M$.
If $p\in N_1\cap N_2$ and $T_pN_1=T_pN_2$, then $N_1=N_2$.
\end{lemma}
Since $f(F^1_b)$, $F^2_b$ are totally geodesic in a complete
manifold, they are complete. By the definition of $f$, $f(p)=q$,
 $f(p)\in f(F^1_b)\cap F^2_b$.  By \eqref{e:step1},
$T_{f(p)}(f(F^1_b))=T_{f(p)}F^2_b$, which, by Lemma
\ref{l:esc-oneill}, implies that $f(F^1_b)=F^2_b$. It follows that
$(\pi\circ f)(z)=\pi_2(z)$ for every $z\in F^1_b$ and that
$T_{f(z)}f(F^1_b)=T_{f(z)}F^2_b$  for every $z\in F^1_b$. Hence,
$f_{*z}(\mathcal H^1_z)=\mathcal H^2_{f(z)}$ for every $z\in F^1_b$.
Since
 $f_{*p}=(\pi_{2*q}|_{\mathcal H^2})^{-1}\circ(\pi_{1*p}|_{\mathcal H^1})$ and since
 every
 vector  of  $\mathcal B^1$ and $\mathcal B^2$  is basic
 along
 $F^1_b$ and $F^2_b$, respectively,
$f_{*z}(A^1_EF)=A^2_{f_{*z}E}f_{*z}F$ for  every $E,F\in T_zH^a_l$
and every $z\in F^1_b$.

{\bf Proof of Step 3.} Let $\gamma:[0,1]\to B$ be a geodesic in $B$
starting from $b=\gamma(0)$. Let $c=\gamma(1)$. For any $z\in
F^1_b$, $w\in F^2_b$ we denote by
 $\gamma^1_z:[0,1]\to H^a_l$ and $\gamma^2_w:[0,1]\to H^a_l$ the $\pi_1$- and $\pi_2$-horizontal lifts of $\gamma$
 starting from  $ z=\gamma^1_{z}(0)$ and from $ w=\gamma^1_{w}(0)$,
 respectively.
 Note that the global existence of the horizontal lifts is ensured by
 the Ehresmann-completeness of $\mathcal H$.
 Let
 $\tau^1_\gamma:F^1_b\to F^1_c$ and $\tau^2_\gamma:F^2_b\to F^2_c$
 be the holonomy diffeomeorphisms
of $\gamma$, given by
 $\tau^1_\gamma(z)=\gamma^1_z(1)$ and $\tau^2_\gamma(w)=\gamma^2_w(1)$, respectively
 (see \cite{gromoll-walschap,bes}).
  A nice fact to point out is that
  $\tau^1_\gamma$ and $\tau^2_\gamma$ are
 isometries since the fibres are totally geodesic \cite{hermann,bes}.
 Now, we prove that
$f\circ\tau^1_\gamma(z)=\tau^2_\gamma\circ f(z)$ for any $z\in
F^1_b$.

The geodesic $f\circ\gamma^1_z$ is $\pi_2$-horizontal if its initial
velocity is (cf. \cite{bes,esc}). We see that
\begin{eqnarray}
\frac{d}{dt}(f\circ\gamma^1_z)\big|_{t=0}=f_{*z}(\dot\gamma^1_z(0))\in
f_{*z}(\mathcal H^1_z)=\mathcal H^2_{f(z)}.
\end{eqnarray}
Thus $\gamma^2_{f(z)}=f\circ\gamma^1_z$ for any $z\in F^1_b$, which
can be reinterpreted as $f\circ\tau^1_\gamma(z)=\tau^2_\gamma\circ
f(z)$.   Therefore, $f(F^1_c)=F^2_c$, hence $f_{*z}(\mathcal
H^1_z)=\mathcal H^2_{f(z)}$ and $\pi_2\circ f(z)=\pi_1(z)$ for any
$z\in F^1_c$.

We now check that $f$ preserves the O'Neill integrability tensors.
Let $X'(t),$ $Y'_1(t),$ $\cdots,$ $Y'_r(t),$ $L'_1(t),$ $\cdots,$
$L'_{k-1}(t)$ be the parallel transports along $\gamma$ of
 $\pi_{1*}X^1,$
 $\pi_{1*}Y^1_1,$ $\cdots,$ $\pi_{1*}Y^1_r,$
 $\pi_{1*}L^1_1,$ $\cdots,$ $\pi_{1*}L^1_{k-1}.$
Let $(X^1(t),$ $Y^1_1(t), \cdots,Y^1_r(t), L^1_1(t),\cdots,
L^1_{k-1}(t))$ and $(X^2(t),$ $Y^2_1(t),
 \cdots,Y^2_r(t),L^2_1(t),\cdots,$ $L^2_{k-1}(t))$
 be the $\pi_1$- and
$\pi_2$-horizontal lifts of $$(X'(t), Y'_1(t), \cdots,Y'_r(t),
L'_1(t),\cdots, L'_{k-1}(t))$$ along $F^1_{\gamma(t)}$ and
$F^2_{\gamma(t)}$, respectively. Set $v_i(t)=A^1_{X^1(t)}Y^1_i(t)$
and
 $w_i(t)=A^2_{X^2(t)}Y^2_i(t)$.
  Fixing $z\in F^1_b$, we simply define   $\gamma^1=\gamma^1_z$.
 We need to establish the following technical lemma.

\begin{lemma}\label{l:parallel-trans}
 {\rm(i)} We have $v^1(\nabla_{\dot\gamma^1(t)}A^1_{X^1(t)}{Y^1_i(t)})=0$ and
        $v^1(\nabla_{\dot\gamma^1(t)}A^1_{L^1_\alpha(t)}{L^1_\beta(t)})=0$, for any $i,\alpha,\beta$.
\begin{itemize}
\item[{\rm(ii)}] The basis
 $\{v_1(t),\cdots,v_r(t)\}$ is an orthonormal basis of vector fields on the
 fibre $F^1_{\gamma(t)}$.
\item[{\rm(iii)}] We have $h^1(\nabla_{\dot\gamma^1(t)}A^1_{L^1_\alpha(t)}{v_i(t)})=0$.
\item[{\rm(iv)}] The vector field
 $\pi_{1*}(A^1_{L^1_{\alpha}(t)}v_i(t))$ is the parallel
 transport of $\pi_{1*}(A^1_{L^1_{\alpha}}v_i)$.
\item[{\rm(v)}] The basis $\mathcal
B^1(t)=\{L^1_0(t),A^1_{L^1_0(t)}v_1(t)\cdots,A^1_{L^1_0(t)}v_r(t),\cdots
L^1_{k-1}(t),A^1_{L^1_{k-1}(t)}v_1(t)\cdots,\\ A^1_{L^1_{k-1}(t)}v_r(t)\}$
is an orthonormal basis of $\mathcal H^1_{\gamma^1(t)}$, and
moreover $A^1_{L^1_\alpha(t)}{L^1_\beta(t)}=0$ for any $\alpha$ and
$\beta$.
\end{itemize}
\end{lemma}
\begin{proof}[Proof of Lemma \ref{l:parallel-trans}]
(i) Since $H^a_l$ has constant curvature, by Proposition
\ref{p:oneill}(b), we get
\begin{eqnarray*}
0&=&R(X^1(t),Y^1_i(t),\!\dot{\ \gamma^1},U)
 =g((\nabla_{\!\dot{\
\gamma^1}}A^1)_{X^1(t)}Y^1_i(t),U)\\
 &=&g(\nabla_{\!\dot{\
\gamma^1}}A^1_{X^1(t)}Y^1_i(t),U)-g(A^1_{\nabla_{\!\dot{\
\gamma^1}}{X^1(t)}}Y^1_i(t),U)-g(A^1_{X^1(t)}\nabla_{\!\dot{\
\gamma^1}}Y^1_i(t),U)\\
&=&g(\nabla_{\!\dot{\ \gamma^1}}A_{X^1(t)}Y^1_i(t),U).
\end{eqnarray*}
Therefore
 $v^1(\nabla_{\dot\gamma^1(t)}A^1_{X^1(t)}{Y^1_i(t)})=0$. Similarly, we
 get $v^1(\nabla_{\dot\gamma^1(t)}A^1_{L^1_\alpha(t)}{L^1_\beta(t)})=0$.

 (ii)  We simply have
\begin{eqnarray}
\!\dot{\ \gamma^1}(t)g(v_i(t),v_j(t))=g(v^1\nabla_{\!\dot{\
\gamma^1}(t)}v_i,v_j)+g(v_i,v^1\nabla_{\!\dot{\ \gamma^1}(t)}v_j)
 =0,
\end{eqnarray}
which implies that $g(v_i(t),v_j(t))$ is constant along
$\gamma^1(t)$ and thus $\{v_i(t)\}_{1\leq i\leq r}$ is an
orthonormal basis.

  (iii)  Using the fact that $(\nabla_{E_1}A)_{E_2}$ is
 skew-symmetric with respect to $g$
 (see \cite{bes}), and that the total space has constant curvature, by Proposition \ref{p:oneill}(b), we have
\begin{eqnarray*}
0&=&R(L^1_\alpha(t),Z,\!\dot{\ \gamma^1},v_i(t))
 =g((\nabla_{\!\dot{\
\gamma^1}}A^1)_{L^1_\alpha(t)}Z,v_i(t))=-g(Z,(\nabla_{\!\dot{\
\gamma^1}}A^1)_{L^1_\alpha(t)}v_i(t))\\
 &=&
 -g(Z,\nabla_{\!\dot{\ \gamma^1}}A^1_{L^1_\alpha(t)}v_i(t))
 +g(Z,A^1_{\nabla_{\!\dot{\ \gamma^1}}{L^1_\alpha(t)}}v_i(t))
 +g(Z,A^1_{L^1_\alpha(t)}v^1\nabla_{\!\dot{\ \gamma^1}}v_i(t))
\\
&=&-g(Z,\nabla_{\!\dot{\ \gamma^1}}A^1_{L^1_\alpha(t)}v_i(t)),
\end{eqnarray*}
which implies  (iii).
 (iii), we simply have
$\nabla'_{\dot\gamma(t)}\pi_{1*}(A^1_{L^1_\alpha(t)}{v_i(t)})=0$.

 (v)  By  (iv), we have that $\mathcal B^1(t)$ is an orthonormal basis
of $\mathcal H^1_{\gamma^1(t)}$. By $(i)$, we get
\begin{eqnarray*}
\!\dot{\
\gamma^1}(t)g(A^1_{L^1_\alpha(t)}{L^1_\beta(t)},v_i(t))=g(v^1\nabla_{\!\dot{\
\gamma^1}(t)}A^1_{L^1_\alpha(t)}{L^1_\beta(t)},v_i)+g(A^1_{L^1_\alpha(t)}{L^1_\beta(t)},v^1\nabla_{\!\dot{\
\gamma^1}(t)}v_i)
 =0,
\end{eqnarray*}
which implies that
$g(A^1_{L^1_\alpha(t)}{L^1_\beta(t)},v_i(t))=g(A^1_{L^1_\alpha(0)}{L^1_\beta(0)},v_i(0))=0$,
for any $i$. Therefore, $A^1_{L^1_\alpha(t)}{L^1_\beta(t)}=0$.
\end{proof}
Similar results hold for $\pi_2$. In particular,
$\pi_{2*}(A^2_{L^2_{\alpha}(t)}w_i(t))$ is the parallel
 transport of $\pi_{2*}(A^2_{L^2_{\alpha}}w_i)$. From Step 2,
 $\pi_{1*}(A^1_{L^1_{\alpha}}v_i)=\pi_{2*}(A^2_{L^2_{\alpha}}w_i)$, and
 therefore their parallel transports must be equal to each other:
\begin{eqnarray}
 \pi_{1*z}(A^1_{L^1_{\alpha}(t)}v_i(t))=\pi_{2*f(z)}(A^2_{L^2_{\alpha}(t)}w_i(t)),
\end{eqnarray}
 and that can be rewritten as
 $f_{*z}(A^1_{L^1_{\alpha}(t)}v_i(t))=A^2_{L^2_{\alpha}(t)}w_i(t) $.
 Using an argument similar to Step 2 for the special bases $\mathcal B^1(t)$
 and $\mathcal B^2(t)$, we simply get $f_{*z}(A^1_EF)=A^2_{f_*E}f_*F$ for
 any $E,F\in\mathcal B^1(t)$.

{\bf Proof of Step 4.} Let $x$ be an arbitrary point in $H^a_l$.
Since $H^a_l$ is connected, there exists a broken geodesic
$\gamma(t)$ in $B$ connecting $b$ and $\pi_1(x)$ (see
\cite[p.~72]{onei}). Applying successively Step 3 to each smooth
piece of the broken geodesic, we see that   $(\star)$ is satisfied
at every point $z\in F_{\gamma(t)}$, for every $t$; in particular,
$(\star)$ holds at $x$.
\end{proof}

\begin{remark}
A very important result due to Escobales is the criterion of
equivalence of two Riemannian submersions, which states that if
$\pi_1,\pi_2:M\to B$ are Riemannian submersions with connected
totally geodesic fibres from a connected complete Riemannian
manifold onto a Riemannian manifold, and if, for some isometry
$f:M\to M$ the condition $(\star)$ holds at a given point $p\in M$,
then there exists an isometry $\tilde f:B\to B$ such that
$\pi_2\circ f=\tilde f\circ\pi_1$.
 Although the proof of Lemma
\ref{l:parallel-trans}(i) invokes $R(X,Y,Z,U)=0$, a usual hypothesis
in the geometry of transversally symmetric (pseudo-)Riemannian
foliations (see \cite{tondeur}), the proof of Theorem
\ref{t:criterium} relies on the construction of a special basis,
which is specific to a pseudo-Riemannian submersion with totally
geodesic fibres of a non-flat real space form. In Theorems
\ref{t:CB} and \ref{t:AB}, we shall see that Theorems
\ref{t:criterium} can be adapted to the case of pseudo-Riemannian
submersions with (para-)complex, connected, totally geodesic fibres
from a (para-)complex pseudo-hyperbolic space.
\end{remark}

\section{Applications of the main theorem}\label{s6}
\noindent We summarize the results proved in the previous sections.
\begin{proof}[Proof of Theorem \ref{t:main}]
By Theorems \ref{t:constant0}, \ref{t:constant3}, \ref{t:constant2}
and Corollary
 \ref{t:nocayley}, $B$ is isometric to one of the following spaces
$H^8_{4}(-4), H^8(-4),  H^{8}_{8}(-4), \mathbb CH^m_t,\mathbb
AP^m,\mathbb HH^m_t, \mathbb BP^m,$ denoted simply by $B'$. There
exists an isometry $\tilde f:B\to B'$. Let $\pi':M'\to B'$ be the Hopf
pseudo-Riemannian submersion with the base space $B'$ and with $M'$
a pseudo-hyperbolic space. Also, by Theorems \ref{t:constant0},
\ref{t:constant3}, \ref{t:constant2},
  we see that $a=\dim(M')$, $l=\mathrm{index}(M')$, and thus $M'=H^a_l$. By
 Theorem \ref{t:criterium}, $\pi':H^a_l\to B'$ is equivalent to
 $\tilde f\circ\pi:H^a_l\to B'$, namely there exists an isometry $f:H^a_l\to H^a_l$ such that
 $\pi\circ f=\tilde f\circ\pi$. Therefore, $\pi$ and $\pi'$
 are equivalent.
\end{proof}

As a consequence of Theorem \ref{t:main}, we now obtain
classification results for pseudo-Riemannian submersions with
totally geodesic fibres from
 $(a)$ $\mathbb CH^m_t,$   $(b)$ $\mathbb  HH^m_t,$
  $(c)$ $\mathbb AP^m,$
  $(d)$ $\mathbb  BP^m.$
First, we define the following Hopf pseudo-Riemannian submersions
with totally geodesic fibres:
\begin{itemize}
\item[$(a)$] $\pi_{\mathbb C,\mathbb H}:\mathbb
CH^{2m+1}_{2t+1}=H^{4m+3}_{4t+3}/H^1_1\to
 \mathbb HH^{m}_{t}=H^{4m+3}_{4t+3}/H^3_3, $
   given by
$\pi_{\mathbb C,\mathbb H}([zH^1_1])=[zH^3_3]$;
\item[$(b)$] $\pi_{\mathbb C,\mathbb B}:\mathbb
CH^{2m+1}_{m}=H^{4m+3}_{2m+1}/H^1_1\to\mathbb
BP^{m}=H^{4m+3}_{2m+1}/H^3_1,$ given by $\pi_{\mathbb C,\mathbb
B}([zH^1_1])=[zH^3_1]$;
\item[$(c)$] $\pi_{\mathbb A,\mathbb B}:\mathbb AP^{2m+1}=H^{4m+3}_{2m+1}/H^1\to\mathbb
BP^{m}=H^{4m+3}_{2m+1}/H^3_1,$ given by $\pi_{\mathbb A,\mathbb
B}([zH^1])=[zH^3_1]$.
\end{itemize}
The fibres of $\pi_{\mathbb C,\mathbb H}$, $\pi_{\mathbb C,\mathbb
B}$, $\pi_{\mathbb A,\mathbb B}$ are isometric to $\mathbb CH^1_1$,
 $\mathbb CH^1$, $\mathbb AP^1,$ respectively.

\begin{theorem}\label{t:CB}
 If $\pi: \mathbb CH^a_b\to B$ is a pseudo-Riemannian submersion
with connected totally geodesic fibres from a complex
pseudo-hyperbolic space onto a pseudo-Riemannian manifold and if the
fibres are complex submanifolds then $\pi$ is equivalent to one of
the following Hopf pseudo-Riemannian submersions:
$$
(a)\ \pi_{\mathbb C,\mathbb H}:\mathbb CH^{2m+1}_{2t+1}\to\mathbb
HH^{m}_{t} \ \ \ \ \ \
 (b)\ \pi_{\mathbb C,\mathbb B}:\mathbb
CH^{2m+1}_{m}\to\mathbb BP^{m}
$$
\end{theorem}
\begin{proof}
Let $\theta:H^{2a+1}_{2b+1}\to\mathbb CH^a_b$ be the Hopf
pseudo-Riemannian submersion over $\mathbb CH^a_b$. Now, $\pi$ and
$\theta$ are pseudo-Riemannian submersions with totally geodesic
fibres, and by Escobales \cite[Theorem 2.5]{esco} so is $\pi\circ\theta$, to which
we can apply Theorem \ref{t:main}. By our usual assumption
$\dim\mathbb CH^a_b>\dim B$, we see that the dimension of the fibres
of $\pi\circ\theta$ is greater than $1$. Therefore, $\pi\circ\theta$
is equivalent to the Hopf pseudo-Riemannian submersions (c), (d),
(e), (f), (g) of Theorem \ref{t:main}, which implies that $\pi$ must
be of the following forms:
\begin{itemize}
\item[$(i)$]\ $\mathbb CH^{2m+1}_{2t+1}\to\mathbb HH^{m}_{t}$,\ \ \ \ \
 $(ii)$\ $\mathbb CH^{2m+1}_{m}\to\mathbb BP^{m},$
\item[$(iii)$]\ $\mathbb CH^{7}_{3}\to H^8(-4)$, \ \ \ \ \
  $(iv)$\ $\mathbb CH^{7}_{3}\to H^8_{4}(-4)$, \ \ \ \ \
   $(v)$\ $\mathbb CH^{7}_{7}\to H^8_{8}(-4)$,
\end{itemize}
By Nagy \cite[Proposition 4.2]{nagy}, the dimension of the fibres must be
$2$, thus, $(iii)$-$(v)$ are not possible. We refer the reader to
\cite{ran} for a different proof of the non-existence of $(v)$, and
to \cite{dga} for that of $(iii)$.
 Let $\pi_1,\pi_2:
\mathbb CH^{2m+1}_{2t+1}\to\mathbb HH^{m}_{s}$ be two
pseudo-Riemannian submersions with totally geodesic fibres. By
Theorem \ref{t:criterium}, $\pi_1\circ\theta$ and $\pi_2\circ\theta$
are equivalent, and, by the proof of Theorem \ref{t:criterium}, there
exists an isometry $f:H^{4m+3}_{4t+3}\to H^{4m+3}_{4t+3}$ depending
on the choice of  an orthonormal basis $\{v_{1p},v_{2p},v_{3p}\}$ of
$\mathcal V^1_p=\mathrm{Ker}(\pi_1\circ\theta)$, $p\in
H^{4m+3}_{4s+3}$, such that
\begin{eqnarray}\label{e:theta}
\pi_2\circ\theta\circ f=\pi_1\circ\theta.
\end{eqnarray}
 If we choose this orthonormal basis such that $v_{3p}$ is
$\theta$-vertical, then, by a similar argument to the proof of
Theorem \ref{t:criterium}, we see that $f$ sends any $\theta$-fibre
into a $\theta$-fibre, and thus there exists an isometry $\tilde
f:\mathbb CH^{2m+1}_{2s+1}\to \mathbb CH^{2m+1}_{2s+1}$ such $\tilde
f\circ\theta=\theta\circ f$. By \eqref{e:theta}, we get $\pi_2\circ
f=\pi_1$.

A similar argument can be used to show the equivalence of   two
pseudo-Riemannian submersions   $\pi_1,\pi_2:\mathbb
CH^{2m+1}_{m}\to\mathbb BP^{m}$.
\end{proof}

\begin{theorem}\label{t:AB}
If $\pi: \mathbb AP^a\to B$ is a pseudo-Riemannian submersion with
connected totally geodesic fibres from a para-complex projective
space onto a pseudo-Riemannian manifold and if the fibres are
para-complex submanifolds then $\pi$ is equivalent to the Hopf
pseudo-Riemannian submersions:
$$
\ \pi_{\mathbb A,\mathbb B}:\mathbb AP^{2m+1}\to\mathbb BP^{m}.
$$
\end{theorem}
\begin{proof}
Let $\pi_{\mathbb A}:H^{2a+1}_{a}\to \mathbb AP^{a}$ be the Hopf
pseudo-Riemannian submersion over $\mathbb AP^{a}$. One can show by
an analogous argument to \cite[Proposition 4.2]{nagy} that  in the
para-case the fibres are also of dimension $2$. Applying Theorem
\ref{t:main} to $\pi\circ\pi_{\mathbb A}$, we obtain that $\pi$
should be of the form
$$
(i)\ \mathbb AP^{2m+1}\to\mathbb BP^{m}, \text{ or } (ii)\ \mathbb
AP^{4m+3}\to\mathbb HH^{2m+1}_m.
$$
Since the signatures of $\mathbb HH^{2m+1}_m$ and $\mathbb
AP^{4m+3}$ are $(4m+4,4m)$ and  $(4m+3,4m+3)$, respectively, $(ii)$
is not possible. The uniqueness of $(i)$ follows analogously to the
proof of Theorem \ref{t:CB}.
\end{proof}

\begin{remark}
The two twistor spaces $\pi:(Z^\varepsilon,g) \to\mathbb BP^m$ ,
$\varepsilon=\pm  1$ (\cite{alek})
  of the  para-quaternionic K\"ahler manifold $\mathbb BP^n$
 are equivalent to the Hopf pseudo-Riemannian submersions\\
$\pi_{\mathbb C,\mathbb B}:\mathbb CH^{2m+1}_{m} \to\mathbb BP^{m}$
(when $\varepsilon=-1$) and
 $\pi_{\mathbb A,\mathbb B}:\mathbb AP^{2m+1}\to\mathbb
 BP^{m}$ (when $\varepsilon=1$). Here $g$ is the canonical K\" ahler-Einstein (when $\varepsilon=-1$)
or para-K\"ahler-Einstein (when $\varepsilon=1$)
 metric of  $Z^\varepsilon$ (see \cite{alek}). By Alekseevsky and Cort\' es \cite[Theorem 3]{alek},
  there are two Einstein metrics
 in the canonical variation on $Z^\varepsilon$ and only one of them is
 $\varepsilon$-K\" ahler-Einstein.
  Another nice fact is
that the twistor space $\pi:Z\to\mathbb HH^m_t$ of the quaternionic
 K\"ahler manifold $\mathbb HH^m_t$ is equivalent to
$\pi_{\mathbb C,\mathbb H}:
  \mathbb CH^{2m+1}_{2t+1}\to \mathbb HH^{m}_{t}. $
\end{remark}
\begin{corollary}
{\rm (i)} There are no pseudo-Riemannian submersions $\pi:\mathbb
HH^m_t\to B$ with connected quaternionic fibres.

{\rm (ii)}  There are no pseudo-Riemannian submersions $\pi:\mathbb
BP^m\to B$ with connected para-quaternionic fibres.
\end{corollary}

\begin{proof}
 First, we recall that any (para-)quaternionic submanifold of a (para-)quaternionic manifold
 is totally geodesic \cite{alekseevsky-marchiafava}.

 (i)  To obtain a contradiction, suppose that such a submersion $\pi$ exists.
 Let $\pi_{\mathbb H}:H^{4m+3}_{4t+3}\to \mathbb
HH^m_t$ be the Hopf pseudo-Riemannian submersion over $\mathbb
HH^m_t$. By Theorem \ref{t:main},
 $\pi\circ\pi_{\mathbb H}$ is equivalent to one of the following:
 $H^{15}_{7}\to H^8(-4)$, $H^{15}_{7}\to H^8_{4}(-4)$, or
 $H^{15}_{15} \to H^{8}_{8}(-4) $, thus $\pi$ must be of the form
\begin{eqnarray}\label{e:quat}
 (a)\ \ \mathbb HH^3_1\to H^8(-4), \ \  (b)\ \mathbb HH^3_1\to
H^8_4(-4)\
 \text{ or }
 (c)\ \mathbb HH^3_3\to H^8_8(-4).
\end{eqnarray}
  We conclude that the fibres are four-dimensional and that  $\pi\circ\pi_{\mathbb
 C,\mathbb H}:\mathbb CH^{7}_{2t+1}\to H^8_{s}(-4)$, $(t,s)\in\{(1,0),(1,4),(3,8)\}$ are pseudo-Riemannian
 submersions with complex, totally geodesic, six-dimensional fibres, which contradicts Theorem
 \ref{t:CB}.

 The proof of (ii) is analogous to (i).
\end{proof}

\begin{remark} The Ucci topological proof \cite{ucc} of the non-existence of
 (\ref{e:quat}(c))  cannot be extended to (\ref{e:quat}(a)) and (\ref{e:quat}(b)),
  because  $\mathbb HH^3_1$, $H^8(-4)$, $H^8_4(-4)$ have
the homotopy types of $S^4$, a point and $S^4$, respectively.
\end{remark}
\begin{remark}
Unlike the Riemannian submersions from spheres, the
pseudo-Riemannian ones from pseudo-hyperbolic spaces feature less
rigidity when we drop the condition of totally geodesic fibres.
 Particularly, while any Riemannian submersion from a sphere is equivalent to a Hopf one
\cite{wilk}, this is no longer true for the pseudo-Riemannian
submersions from pseudo-hyperbolic spaces. Indeed (cf.
\cite{barros-fernandez}) any pseudo-hyperbolic space $H^a_l$ can
simply be written as a warped product
$H^a_l=(H^{a-l}\times_{f}S^l,g_{H^a_l}),$  via the identification
$\phi:H^{a-l}\times S^l\to H^a_l$, given by
 $\phi((x_0, x), u) = (x_0u, x)$, for every
 $u\in S^l$, $(x_0, x)\in H^{a-l}$, $x_0\in\mathbb R_+,$
 $x\in\mathbb R^{a-l}$. Here $f:H^{a-l}\to \mathbb R_+$ is given
 by $f(x_0,(x_1,\cdots,x_{a-l}))=x_0$, and the metric of the warped
 product is
 $ g_{H^{a-l}}-f^2g_{S^l}$.
Now, the
projection
$$\pi:H^a_l=H^{a-l}\times_{f}S^l\to H^{a-l}$$
 is a pseudo-Riemannian
submersion (with totally umbilical fibres \cite{bes}), which is not
equivalent to a Hopf one, except possibly when
$(a,l)\in\{(3,1),(7,3),(15,7)\}$. The classification problem  of
pseudo-Riemannian submersions from pseudo-hyperbolic spaces remains
open.
\end{remark}
\bibliographystyle{amsplain}

\end{document}